\documentclass[12pt,twoside]{article}

\usepackage{graphics}
\usepackage{latexsym}
\usepackage{amsmath}
\usepackage{amssymb}

\makeatletter
\newtheorem{theorem}{Theorem}[section]
\newtheorem{lemma}[theorem]{Lemma}

\newtheorem{remark}[theorem]{Remark}
\newtheorem{cor}[theorem]{Corollary}
\newtheorem{crit}[theorem]{Criterion}
\newtheorem{prop}[theorem]{Proposition}
\newtheorem{defin}[theorem]{Definition}
\newtheorem{notation}[theorem]{Notation}
\newtheorem{ex}[theorem]{Example}

\newtheorem{constr}[theorem]{Construction}

\def\ps@headings{
 \def\@oddhead{\footnotesize\rm\hfill\runningheadodd\hfill\thepage}
 \def\@evenhead{\footnotesize\rm\thepage\hfill\runningheadeven\hfill}
 \def\@oddfoot{}
 \def\@evenfoot{\@oddfoot}
}

\newcommand{\Prf}{\noindent{\bf Proof}.\quad }

\newcommand{\qed}{\hfill$\Box$}

\makeatother

\def\runningheadeven{On the Ubiquity and Utility of Cyclic Schemes}
\def\runningheadodd{M. Abreu, M.J. Funk, D. Labbate, V. Napolitano}
\title{On the Ubiquity and Utility of Cyclic Schemes}
\author{{\rm M. Abreu,${}^\dagger$\thanks{This research was carried out within the activity
 of INdAM-GNSAGA and supported by the Italian Ministry MIUR  $\;-\;$  $^c$ Correponding author}   $\;$
 M.J. Funk,$^{\dagger \ast c } $$\;$D. Labbate,$^{\ddagger \ast}$ $\;$V. Napolitano$^{\S \ast }$}\\
\\ ${}^\dagger$ \small Dipartimento di Matematica e Informatica -
Universit\`a della
        Basilicata \\
      \small 85100 Potenza -
         Italy (marien.abreu@unibas.it; martin.funk@unibas.it)
          \\
\\ $^\ddagger$ \small Dipartimento di Matematica - Politecnico di Bari
\\ \small Via E. Orabona, 4 - 70125 Bari -
         Italy (labbate@poliba.it)\ \\
         \\ $^\S$ \small Dipartimento di Ingegneria Civile - Seconda Universit\`a di Napoli
\\ \small  via Roma, 29 - 81031 Aversa -
         Italy (vito.napolitano@unina2.it)
         }

      \date{}

\begin{document}
\maketitle
\pagestyle{headings}

\begin{abstract}
Let $k,l,m,n$, and $\mu$ be positive integers. A
$\mathbb{Z}_\mu$--{\it scheme of valency} $(k,l)$ and {\it order} $(m,n)$ is a $m \times
n$ array $(S_{ij})$ of subsets $S_{ij} \subseteq \mathbb{Z}_\mu$
such that for each row and column one has $\sum_{j=1}^n |S_{ij}| =
k $ and $\sum_{i=1}^m |S_{ij}| = l$, respectively. Any such scheme is an algebraic equivalent of
a $(k,l)$--semi--regular bipartite voltage
graph with $n$ and $m$ vertices in the bipartition sets and
voltages coming from the
cyclic group $\mathbb{Z}_\mu$. We are
interested in the subclass of $\mathbb{Z}_\mu$--schemes that are
characterized by the property $a - b + c - d\; \not \equiv \;0$
(mod $\mu$) for all $a \in S_{ij}$, $b \in S_{ih}$, $c \in
S_{gh}$, and $d \in S_{gj}$ where $i,g \in \{1,\ldots,m\}$ and
$j,h \in \{1,\ldots,n\}$ need not be distinct. These
$\mathbb{Z}_\mu$--schemes can be used to represent adjacency
matrices of regular graphs of girth $\ge 5$ and semi--regular
bipartite graphs of girth $\ge 6$. For suitable $\rho, \sigma \in
\mathbb{N}$ with $\rho k = \sigma l$, they also represent
incidence matrices for polycyclic $(\rho \mu_k, \sigma \mu_l)$
configurations and, in particular, for all known Desarguesian
elliptic semiplanes.  Partial projective closures yield {\it mixed
$\mathbb{Z}_\mu$--schemes}, which allow new constructions for Kr\v
cadinac's sporadic  configuration of type $(34_6)$ and Balbuena's
bipartite $(q-1)$--regular graphs of girth $6$ on as few as
$2(q^2-q-2)$ vertices, with  $q$ ranging over prime powers.
Besides some new results, this survey essentially furnishes new
proofs in terms of (mixed) $\mathbb{Z}_\mu$--schemes for ad--hoc
constructions used thus far.
\medskip

Keywords:  (mixed) cyclic schemes, cyclic voltage graphs, (polycyclic) configurations, elliptic semiplanes, small regular graphs with girths $5$ and $6$

\medskip
AMS MSC:  05B20,05B30, 05C50, 05C75, 51A45
\end{abstract}

\section{$\mathbb{Z}_\mu$--Schemes and Cyclic Voltage Graphs}

{\bf Preliminary note.} This paper deals with constructions in
some classes of $(0,1)$--matrices, which turn up as incidence
matrices of configurations or  adjacency matrices of graphs. Even
if basic notions and notations seem to be generally known and
widely used, misunderstandings can arise since precise formal
definitions vary slightly from author to author. So one might be
tempted to fix every notion to the least detail, at the risk of
distracting the reader's attention from the essentially new
concepts. To overcome this dilemma, the reader will find a
synopsis on $(0,1)$-matrices, graphs, and configurations in
Section \ref{verylast}. Notions defined in the synopsis are set up
in italics at their very first appearance in the paper.

\begin{defin}\label{S}  A $\mathbb{Z}_\mu$--{\bf scheme of
order} $(m,n)$ is an $m \times n$ array $M^{(\mu)} = (S_{ij})$ of
subsets $S_{ij} \subseteq \mathbb{Z}_\mu$. The
$\mathbb{Z}_\mu$--scheme $M^{(\mu)}$ has {\bf valency} $(k,l)$ if,
for each row and column, the sums of the cardinalities of the
entries have constant values $k$ and $l$, i.e.
$$\sum_{j=1}^n |S_{ij}| = k \quad \makebox{and} \quad
\sum_{i=1}^m |S_{ij}| = l\,.$$ If each entry has cardinality $\le
1$ and precisely $1$, the scheme is called {\bf simple} and {\bf
full}, respectively. If $m = n$ and $k = l$, we say that
$M^{(\mu)}$ has {\bf order} $n$ and {\bf valency} $k$,
respectively. A  $\mathbb{Z}_\mu$--scheme of order $n$ is said to
be {\bf skew--symmetric} if $S_{ij} = -S_{ji}$ for all $1 \le i,j
\le n$.
\end{defin}

\begin{notation} When writing down a $\mathbb{Z}_\mu$--scheme $M^{(\mu)}$,
the curly braces of the entries will always be omitted.
Accordingly, the empty set $\emptyset = \{\}$ becomes a blank
entry. If necessary, $\mu$ will be mentioned as superscript
$^{(\mu)}$.
\end{notation}
\medskip

A {\it circulant} $(0,1)$--matrix, say $\overline C$, is uniquely
determined by the positions of the entries $1$ in its first row.
This gives rise to a bijective mapping, say $\iota$, from the
class of circulant $(0,1)$--matrices of order $\mu$ onto the power
set of $\mathbb{Z}_\mu$, namely
$$  \overline C = \left(\begin{smallmatrix}
c_0&& c_1&& \hdots&&c_{\mu-2} &&c_{\mu-1}\\\\
c_{\mu-1} && c_0 && c_1 && && c_{\mu-2}\\
\vdots&&c_{\mu-1} && c_0 &&\ddots&&\vdots\\
c_2&& &&\ddots&& \ddots&&c_1\\\\
c_1&&c_2&&\hdots&&c_{\mu-1}&&c_0\end{smallmatrix}\right)\quad
\longmapsto \quad C := \{i \in \mathbb{Z}_\mu\;|\; c_i = 1\,\}\,,
$$ where the empty
set becomes the image of the zero matrix of order $\mu$. When
speaking of positions $(i,j)$ in circulant matrices of order
$\mu$, the indices range over $\{0, \ldots, \mu-1\} =
\mathbb{Z}_\mu$. For later use, the rule determining the inverse
mapping is worthwhile to be stated explicitly:

\begin{lemma} $\;\overline C$ has entry $1$ in position $(i,j)$
if and only if $j - i$ (mod $\mu$) belongs to $C$.\qed\end{lemma}

The mapping induces the following notation for $(0,1)$--block
matrices with circulant blocks (\cite{AFLN}):

\begin{defin}\label{blowup} The {\bf blow--up} of a
$\mathbb{Z}_\mu$--scheme $M^{(\mu)} = (S_{ij})$ is the block
$(0,1)$--matrix $\overline{M^{(\mu)}}$ with square blocks of order
$\mu$ which is obtained from $M^{(\mu)}$ by substituting the
circulant $(0,1)$--matrices $\overline{S_{ij}}$ for the entries
$S_{ij}$.
\end{defin}

In the sequel, the position of an entry in the blow--up
$\overline{M^{(\mu)}}$ will be given in terms of the position
$(i,j)$ of the block $S_{ij}$ and the {\bf local} position
$(i',j')$ within the circulant block $S_{ij}$.

{\it Adjacency matrices} of {\it graphs} are symmetric
$(0,1)$--matrices with entries $0$ on the main diagonal. These two
properties can easily be translated into the language of
$\mathbb{Z}_\mu$--schemes.

\begin{prop} The blow--up $\overline{M^{(\mu)}}$ of a square
$\mathbb{Z}_\mu$--scheme $M^{(\mu)} = (S_{ij})$ of order $n$ is
symmetric if and only if $M^{(\mu)}$ is skew--symmetric.
\end{prop}

\Prf Let $a$ be an element in $S_{ij}$. Then an entry $1$ turns up
in local position $(i',j')$ in the circulant matrix
$\overline{S_{ij}}$ if and only if $j'-i' \equiv a$ (mod $\mu$).
Symmetrically, $1$ appears in local position $(j',i')$ in the
circulant matrix $\overline{S_{ji}}$ if and only if $i'-j' \equiv
-a \in S_{ji}$, numbers taken modulo $\mu$. \qed
\medskip

\begin{cor} The blow--up $\overline{M^{(\mu)}}$ of a square
$\mathbb{Z}_\mu$--scheme $M^{(\mu)} = (S_{ij})$ of order $n$ has
entries $0$ on its main diagonal if and only if $0 \not \in
S_{ii}$ for all $i = 1, \ldots, n$.
\end{cor}

\Prf Apply the above Proof in the case $i = j$ and $i' = j'$ to
see that $1$ is an entry on the main diagonal of
$\overline{S_{ii}}$ if and only if $0 \in S_{ii}$. \qed
\medskip

In the light of these two statements, we call a skew--symmetric
$\mathbb{Z}_\mu$--scheme $M^{(\mu)} = (S_{ij})$ of order $n$ {\bf
admissible} if $0 \not \in S_{ii}$ for all $i = 1, \ldots, n$.
{\it Cyclic voltage graphs} and {\it admissible cyclic voltage
assignments} are surveyed in the beginning of Section
\ref{verylast}.

\begin{remark} \label{1} $(i)$ Any admissible $\mathbb{Z}_\mu$--scheme $M^{(\mu)} = (S_{ij})$ of order $n$
arises from, and gives rise to, a  cyclic voltage graph
$(K,\alpha)$ on $n$ vertices with an
admissible cyclic voltage assignment $\alpha$. Labelling the vertices of $K$  by \/ $1,\ldots,n$,  the rules
$$S_{ij}\, := \,\biggm\{a \in \mathbb{Z}_\mu\,\biggm |\,\alpha(e)
=
\begin{array}{r} a\\-a\\\end{array} \makebox{\rm for some edge}\;
e \in EK\;\makebox{\rm running from} \,
\begin{array}{l}i\;
\makebox{\rm to}\;j \\
j\; \makebox{\rm to}\;i \\\end{array} \biggm\}
$$  for $i \ne j$ as well as $$S_{ii}\; := \;\{a , -a \in
\mathbb{Z}_\mu\;|\; \alpha(e) = a\;\makebox{\rm for
an}\;i\makebox{\rm --based loop} \; e \in EK\; \}
$$ construct $M^{(\mu)}$  from $(K,\alpha)$. Vice versa, given
$M^{(\mu)} = (S_{ij})$, let $K$ be the general graph with vertex
set $VK := \{1, \ldots,n\}$ where $|S_{ij}|$ edges run from $i$ to
$j$ with distinct voltages $a \in  S_{ij}$ and eventually a vertex
$i$ is base of $\frac12 |S_{ii}|$ loops with voltage $\pm b \in
S_{ii}$. Both constructions do comply with admissibility.

\noindent $(ii)$ An arbitrary $\mathbb{Z}_\mu$--scheme $M^{(\mu)}
= (S_{ij})$ of order $(m,n)$ arises from, and gives rise to, a
bipartite cyclic voltage graph on $m$ white and $n$ black
vertices, with an admissible cyclic voltage assignment. Denote by
$- M^{(\mu)}$ the $\mathbb{Z}_\mu$--scheme obtained from
$M^{(\mu)}$ by substituting each entry with its opposite element
in $\mathbb{Z}_\mu$, and let $O_\nu$ be the trivial
$\mathbb{Z}_\mu$--scheme of order $\nu$ all of whose entries are
$\emptyset$. Then
$$\begin{pmatrix} O_m & M^{(\mu)}\\(- M^{(\mu)})^T & O_n \end{pmatrix}$$
is an admissible $\mathbb{Z}_\mu$--scheme of order $m+n$ and $(i)$
applies, both constructions being compatible with bipartite
(general) graphs.
\end{remark}

\begin{prop} If $M^{(\mu)}
= (S_{ij})$  is an admissible $\mathbb{Z}_\mu$--scheme with
associated voltage graph $(K,\alpha)$, the blow--up
$\overline{M^{(\mu)}}$ is an adjacency matrix of the lift of $K$
through $\mathbb{Z}_\mu$ via $\alpha$.
\end{prop}

\Prf Order the vertices $(i,a) \in VK \times \mathbb{Z}_\mu$
lexicographically with respect to the natural orders $1, \ldots,
n$ for the vertices in $K$ and $0,1, \ldots, \mu-1$ for the
elements in $\mathbb{Z}_\mu$. \qed \medskip

Involving some regularity condition, Remark \ref{1} reads:

\begin{prop} A $\mathbb{Z}_\mu$--scheme  of order $(m,n)$
and valency $(k,l)$ is equivalent to a $(k,l)$--semi--regular bipartite
voltage graph on $n$ white and $m$ black vertices with voltages
from the cyclic group $\mathbb{Z}_\mu$, while an admissible
$\mathbb{Z}_\mu$--scheme  of order $n$ and valency $k$ is
equivalent to a $k$--regular cyclic voltage graph on $n$
vertices with an admissible voltage assignment. \qed
\end{prop}

\begin{ex}{\rm The celebrated
 Petersen graph can be seen as a lift of the dumbbell graph through
$\mathbb{Z}_5$:
$$
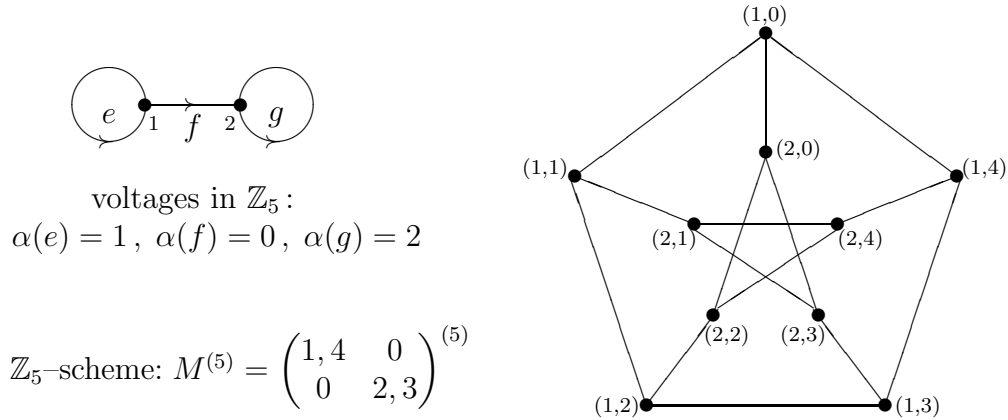
\begin{picture}(20,9)
\put(11,6){\circle*{.3}} \put(12.5,1.2){\circle*{.3}}
\put(17.5,1.2){\circle*{.3}} \put(19,6){\circle*{.3}}
\put(15,9){\circle*{.3}} \put(13.9,3.1){\circle*{.3}}
\put(16.1,3.1){\circle*{.3}} \put(13.5,5){\circle*{.3}}
\put(16.5,5){\circle*{.3}} \put(15,6.5){\circle*{.3}}

\put(15,9.3){\makebox(0,0){$^{(1,0)}$}}
\put(19.6,6.1){\makebox(0,0){$^{(1,4)}$}}
\put(18.2,1.1){\makebox(0,0){$^{(1,3)}$}}
\put(11.9,1.1){\makebox(0,0){$^{(1,2)}$}}
\put(10.4,6.1){\makebox(0,0){$^{(1,1)}$}}
\put(14.2,2.6){\makebox(0,0){$^{(2,2)}$}}
\put(17,4.6){\makebox(0,0){$^{(2,4)}$}}
\put(13.1,4.6){\makebox(0,0){$^{(2,1)}$}}
\put(15.8,2.6){\makebox(0,0){$^{(2,3)}$}}
\put(15.7,6.5){\makebox(0,0){$^{(2,0)}$}}

\put(12.5,1.2){\line(1,0){5}} \put(12.5,1.2){\line(-1,3){1.6}}
\put(12.5,1.2){\line(3,4){1.5}} \put(17.5,1.2){\line(1,3){1.6}}
\put(17.5,1.2){\line(-3,4){1.5}} \put(11,6){\line(4,3){4}}
\put(19,6){\line(-4,3){4}} \put(13.5,5){\line(-5,2){2.5}}
\put(16.5,5){\line(5,2){2.5}} \put(15,6.5){\line(0,1){2.5}}
\put(13.5,5){\line(1,0){3}}

\put(15,6.4){\line(-1,-3){1.1}} \put(15,6.4){\line(1,-3){1.1}}

\put(13.5,4.9){\line(3,-2){2.5}} \put(16.5,4.9){\line(-3,-2){2.5}}

\put(1.25,7.5){\circle{1.5}}\put(4.75,7.5){\circle{1.5}}
\put(2,7.5){\circle*{.3}}\put(4,7.5){\circle*{.3}}
\put(2,7.5){\line(1,0){2}}

\put(2.2,7){\makebox(0,0){$^1$}}\put(3.8,7){\makebox(0,0){$^2$}}
\put(3,7.4){\makebox(0,0){$^\succ$}}\put(1.2,6.625){\makebox(0,0){$^\succ$}}\put(4.7,6.625){\makebox(0,0){$^\succ$}}
\put(3,7){\makebox(0,0){$f$}}\put(1.25,7.25){\makebox(0,0){$e$}}\put(4.75,7.25){\makebox(0,0){$g$}}
\put(3,5.5){\makebox(0,0){voltages in $\mathbb{Z}_5$\,:}}
\put(3.5,4.7){\makebox(0,0){$\alpha(e) = 1\,,\; \alpha(f) =
0\,,\;\alpha(g) = 2$}}
\put(4,2){\makebox(0,0){$\mathbb{Z}_5$--scheme:\;$M^{(5)} = \begin{pmatrix} 1,4 & 0 \\
0 & 2,3\\
\end{pmatrix}^{(5)}$}}
\end{picture}$$

}\end{ex}

\section{$J_2$--Free $\mathbb{Z}_\mu$--Schemes}

As usual, let $J_n$ denote the square matrix all of whose entries
are $1$. Then $J_2$ is the {\it incidence matrix} of a {\bf
di-gon}, i.e. the structure made up by two distinct points $p_1$,
$p_2$, two distinct lines $L_1$, $L_2$, and all four incidences
$p_i | L_j$ with $i,j \in \{1,2\}$. Di-gons are forbidden
substructures of {\it configurations}. Thus, disregarding {\it
regularity} conditions, incidence matrices of configurations
 are $(0,1)$--matrices characterized by the
following property:

\begin{defin} A $(0,1)$--matrix
is called {\bf $J_2$--free} if every $2 \times 2$ submatrix
contains at least one entry $0$. In a figurative sense,  a
$\mathbb{Z}_\mu$--scheme $M^{(\mu)}$ is said to be {\bf
$J_2$-free} if its blow--up $\overline{M^{(\mu)}}$ is so.
\end{defin}

In \cite{AFLN,AFLNA,AFLN5,FLN} such matrices were called \lq \lq
linear."

\begin{crit}\label{criterion 1}
A $\mathbb{Z}_\mu$--scheme $M^{(\mu)} = (S_{ij})$ of order $(m,n)$
is $J_2$-free if and only if for all (not necessarily distinct) $1
\le i,g \le m$ and $1 \le j,h \le n$ and all $a \in S_{ij}$, $b
\in S_{ih}$, $c \in S_{gh}$, and $d \in S_{gj}$ one has
$$(\dagger)\quad\quad a - b + c - d\; \not \equiv
\;0 \;(\makebox{mod}\; \mu)\,.$$
 \end{crit}

\Prf To prove sufficiency, assume that $\overline{M^{(\mu)}}$ has
a sub--matrix $J$ of order $2$ all of whose entries are $1$. By
construction, the  upper left $1$ in $J$ appears as an entry in
local position $(i',j')$ in the block $\overline{S_{ij}}$, for
some $i',j' \in\{0,\ldots,\mu-1\}$, $i \in \{1,\ldots,m\}$, and $j
\in \{1,\ldots,n\}$. This, in turn, implies that $a:\equiv j'-i'$
(mod $\mu$) is an element of the set $S_{ij}$. Analogously,  the
upper right, lower right, and lower left entry $1$ in $J$ arise
from entries $1$ in local positions
$$\begin{array}{cccc}
(i',h')&\makebox{in the block } S_{ih}&
\Longrightarrow & b:\equiv h'-i' \in S_{ih}\,,\\
(g',h')&\makebox{in the block } S_{gh}&
\Longrightarrow & c:\equiv h'-g' \in S_{gh}\,,\\
(g',j')&\makebox{in the block } S_{gj}&
\Longrightarrow & d:\equiv j'-g' \in S_{gj}\,,\\
\end{array}$$ differences taken modulo $\mu$. Subtracting
the second and fourth congruences from the sum of the first and
third, we obtain  $ 0 \equiv a - b + c - d $ (mod $\mu$), a
contradiction.

To prove necessity, suppose that there exist (not necessarily
distinct) $i,g \in \{1,\ldots,m\}$ and $j,h \in \{1,\ldots,n\}$
such that for some $a \in S_{ij}$, $b \in S_{ih}$, $c \in S_{gh}$,
and $d \in S_{gj}$ one has
$$(\ddagger) \qquad a - b + c - d\;  \equiv \;0
\;(\makebox{mod}\; \mu)\,.$$ In the first row of
$\overline{S_{ij}}$ and $\overline{S_{ih}}$, there are entries $1$
in local positions $(0,a)$ and $(0,b)$, respectively. Now consider
the circulant $(0,1)$--block $\overline{S_{gj}}$. Since $d \in
S_{gj}$, there exists a  row index $j' \in \{0, \ldots, \mu-1\}$,
such that $\overline{S_{gj}}$ has an entry $1$ in local position
$(j',a)$, namely $j' :\equiv a - d$ (mod $\mu$). Then $(\ddagger)$
implies $j' \equiv c - b$ (mod $\mu$). Hence $\overline{S_{gh}}$
has an entry $1$ in position $(j',b)$ and $\overline {M^{(\mu)}}$
contains a $2 \times 2$ submatrix all of whose entries are $1$, a
contradiction. \qed
\medskip

Condition $(\dagger)$ has some repercussions on non--empty entries
$S_{ij}$:

\begin{cor}
Each entry $S_{ij} \ne \emptyset$ of a $J_2$--free
$\mathbb{Z}_\mu$--scheme $M^{(\mu)} = (S_{ij})$ is a {\rm
deficient cyclic difference set}.
\end{cor}

\Prf Apply condition $(\dagger)$ in the case that $i = g$ and $j =
h$: all the differences $a - b$, $c - d$ with $a \ne b$ and $c \ne
d$ are distinct in pairs. \qed
\medskip

\begin{cor}
Suppose that the $\mathbb{Z}_\mu$--scheme $M^{(\mu)} = (S_{ij})$
of order $(m,n)$ is $J_2$--free. Then, for all $i,g \in
\{1,\ldots,m\}$ and $j,h \in \{1,\ldots,n\}$, the differences
covered by either $S_{ij}$ and $S_{ih}$ or by $S_{ij}$ and
$S_{gj}$ are pairwise distinct.
\end{cor}

\Prf Apply condition $(\dagger)$ in the case that either $i = g$
or $j = h$. \qed \medskip

\section{Polycyclic Configurations and $\mathbb{Z}_\mu$--Schemes}
\label{polycyclic}

Boben and Pisanski \cite{BoPi} call an $(m_k,n_l)$ configuration
$\cal C$ {\bf polycyclic}  or $\mu$--{\bf cyclic} if $\cal C$
admits a cyclic automorphism of order $\mu$ whose orbits partition
both the point set and the line set of $\cal C$ into subsets of
size $\mu$. This definition makes sense only if $1 < \mu \;|\;
gcd(m,n)$ and  $m = \rho \mu$,  $n = \sigma \mu$ for suitable
$\rho, \sigma \in \mathbb{N}$. {\it Cyclic} configurations $(n_k)$
are $n$--cyclic.

Incidence matrices reveal the polycyclic structure of a
configuration if a suitable labelling matches with the orbits
under the cyclic automorphism.
\medskip

\begin{prop}
A $(\rho \mu_k, \sigma \mu_l)$ configuration $\cal C$ is
polycyclic if and only if it admits an incidence matrix $\overline
{M^{(\mu)}}$  obtained by blowing up a $\mathbb{Z}_\mu$--scheme
$M^{(\mu)}$ of valency $(k,l)$ and order $(\rho,\sigma)$.
\end{prop}

\Prf  Sufficiency is guaranteed by the very construction:
$\overline {M^{(\mu)}}$ admits a cyclic automorphism of order
$\mu$, namely the simultaneous action of the intrinsic cyclic
automorphism on each circulant block of $\overline {M^{(\mu)}}$.
This automorphism induces an automorphism of $\cal C$, whose
orbits partition both the point set and the line set of $\cal C$
into $\rho$  and $\sigma$ subsets of size $\mu$, respectively.
Hence $\cal C$ is polycyclic.

To prove necessity, suppose that $\cal C$ is polycyclic with
respect to some automorphism $\varphi$ of order $\mu$. Choose
representatives $p_0,\hdots,p_{m-1}$ in the point set and
$L_0,\hdots,L_{m-1} $ in the line set of $\cal C$ for the orbits;
i.e. using the abbreviation $(i) := \varphi^i$, one has the
following cycle decompositions:
$$\begin{array}{rcl}(p_0^{(0)}\;p_0^{(1)}\;p_0^{(2)}\,\hdots \,
p_0^{(\mu-1)} )&\hdots &(p_{\rho-1}^{(0)}\;p_{\rho-1}^{(1)}\;
p_{\rho-1}^{(2)}\,\hdots \,p_{\rho-1}^{(\mu-1)} )\\
(L_0^{(0)}\;\,L_0^{(1)}\;\,L_0^{(2)}\,\hdots \;\, L_0^{(\mu-1)}\,
)&\hdots &(L_{\sigma-1}^{(0)}\;\,L_{\sigma-1}^{(1)}\;\,
L_{\sigma-1}^{(2)}\,\hdots \;\,L_{\sigma-1}^{(\mu-1)}\,
)\\\end{array}$$ Let $\overline{M_\varphi}$ be the incidence
matrix for $\cal C$ obtained when labelling its points and lines
according to the above cycle decompositions of $\varphi$.
Interpret
$$\overline{M_\varphi }\quad = \quad
\left(\begin{smallmatrix}\overline{M_{0,0}}& \overline{M_{0,1}}&
\hdots&\overline{M_{0,m-1}} \\\overline{M_{1,0}}&
\overline{M_{1,1}}& \hdots&\overline{M_{1,m-1}}
\\\vdots &\vdots& \ddots&\vdots\\\overline{M_{m-1,1}}&
\overline{M_{m-1,2}}& \hdots&\overline{M_{m-1,m-1}}
\\\end{smallmatrix}\right)$$ as an $m \times m$ block matrix with
square blocks $\overline{M_{i,j}}$ of order $\mu$.

Then each block $\overline{M_{i,j}}$ is a circulant
$(0,1)$--matrix (which might also be a copy of the zero matrix):
in fact, if, for some $i,j,s,t \in \{0, \ldots,\mu-1\}$, the point
$p_i^{(s)}$ and the line $L_j^{(t)}$ are incident, so are their
$\varphi$--images, i.e. $p_i^{(s+1)}$ is incident with
$L_j^{(t+1)}$, apices taken modulo $\mu$; this implies that for
any entry $1$ in position $(s,t)$ in the block $M_{i,j}$, there
exist entries $1$ also in the positions $(s + z,t + z)$ for $z =
1, \ldots, \mu-1$, numbers again taken modulo $\mu$; Hence
$\overline{M_{i,j}}$ is a circulant matrix. A block $M_{i,j}$ has
all its entries $0$ if for each $s \in \{0, \ldots\mu-1\}$ the
point $p_i^{(s)}$ is not incident with any line $L_j^{(t)}$ with
$t \in \{0, \ldots,\mu-1\}$.

$\overline{M_\varphi}$ can be seen as the blow--up of a
$\mathbb{Z}_\mu$--scheme of order $(\rho,\sigma)$, say
$M_\varphi$, obtained by applying $\iota$ to the blocks of
$\overline{M_\varphi}$. The valency of $M_\varphi$ is $(k,l)$, since
$\overline{M_\varphi}$ has exactly $k$ and $l$ entries $1$ in each
row and column, respectively. \qed
\medskip

\begin{ex}{\rm The Cremona--Richmond configuration (\cite{Cox}, represented
geometrically in the figure below) is a $5$--cyclic $(15_3)$
configuration.
$$
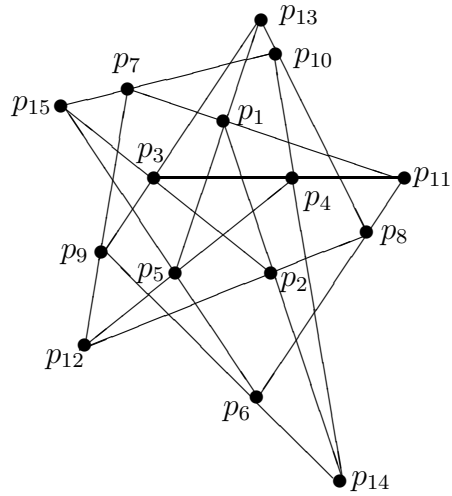
\begin{picture}(20,10)
\put(9,4.5){\circle*{.3}}\put(8.5,4.5){\makebox(0,0){$p_5$}}
\put(9,4.5){\line(1,3){1.8}}
\put(9,4.5){\line(-5,-4){2}}\put(9,4.5){\line(5,4){2.5}}
\put(9,4.5){\line(-2,3){2.3}}\put(9,4.5){\line(2,-3){1.75}}

\put(6.6,8){\circle*{.3}}\put(6,8){\makebox(0,0){$p_{15}$}}
\put(6.6,8){\line(4,1){4.6}}

\put(11,4.5){\circle*{.3}}\put(11.5,4.3){\makebox(0,0){$p_2$}}
\put(11,4.5){\line(-5,4){4.3}}
\put(11,4.5){\line(-5,-2){4}}\put(11,4.5){\line(5,2){2}}

\put(8.55,6.5){\circle*{.3}}\put(8.5,7){\makebox(0,0){$p_3$}}
\put(8.55,6.5){\line(1,0){5.2}}\put(8.55,6.5){\line(2,3){2.2}}
\put(8.55,6.5){\line(-2,-3){1.15}}

\put(10.8,9.8){\circle*{.3}}\put(11.6,9.9){\makebox(0,0){$p_{13}$}}
\put(10.8,9.8){\line(1,-2){2.25}}

\put(11.45,6.5){\circle*{.3}}\put(12,6){\makebox(0,0){$p_4$}}
\put(11.5,6.5){\line(-1,6){.45}}\put(11.45,6.5){\line(1,-6){1.05}}

\put(10,7.7){\circle*{.3}}\put(10.6,7.8){\makebox(0,0){$p_1$}}
\put(10,7.7){\line(1,-3){2.5}}\put(10,7.7){\line(3,-1){3.7}}
\put(10,7.7){\line(-3,1){2.1}}

\put(13.8,6.5){\circle*{.3}}\put(14.4,6.5){\makebox(0,0){$p_{11}$}}
\put(13.8,6.5){\line(-2,-3){3.1}}

\put(12.45,0.15){\circle*{.3}}\put(13.1,0.15){\makebox(0,0){$p_{14}$}}
\put(12.3,0.2){\line(-1,1){4.8}}

\put(7.1,3){\circle*{.3}}\put(6.7,2.7){\makebox(0,0){$p_{12}$}}
\put(7.1,2.9){\line(1,6){.9}}

\put(7.45,4.95){\circle*{.3}}\put(6.9,4.9){\makebox(0,0){$p_9$}}
\put(8,8.35){\circle*{.3}}\put(8,8.9){\makebox(0,0){$p_7$}}
\put(11.1,9.1){\circle*{.3}}\put(11.9,9){\makebox(0,0){$p_{10}$}}
\put(10.7,1.9){\circle*{.3}}\put(10.3,1.6){\makebox(0,0){$p_6$}}
\put(13,5.35){\circle*{.3}}\put(13.6,5.3){\makebox(0,0){$p_8$}}
\end{picture}$$
The permutation
$$(1\;2\;3\;4\;5)(6\;7\;8\;9\;10)(11\;12\;13\;14\;15)\,,$$
acting on the indices of the points, induces an automorphism of
order $5$, which partitions both the point and line sets into
three orbits of length $5$ each. Choose the points $p_1, p_6,
p_{11}$ and the lines
$\{p_3,p_4,p_{11}\}$,$\{p_7,p_{10},p_{15}\}$, $\{p_1,p_7,p_{11}\}$
as first elements in each orbit. Then the resulting incidence
matrix is the blow--up of the $\mathbb{Z}_5$--scheme

$$M_{CR}^{(5)}\quad := \quad \left(\begin{smallmatrix}
2,3 &&& &
0\\
&&1,4 && 4\\0&&1&&0\\\end{smallmatrix}\right)^{(5)}.$$ Note that
the associated bipartite cyclic voltage graph differs slightly
from the one given in \cite[Figure 6]{BoPi} or in \cite[Figure
4(b)]{PBMOG}:
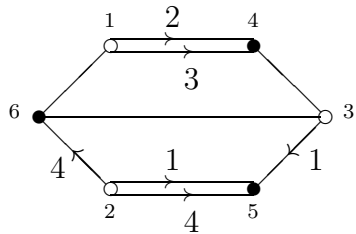
$$\begin{picture}(20,5)
\put(6,2.5){\circle*{.3}}
\put(7.5,4){\circle{.3}}\put(7.5,1){\circle{.3}}
\put(10.5,4){\circle*{.3}} \put(10.5,1){\circle*{.3}}
\put(12,2.5){\circle{.3}}
\put(6,2.5){\line(1,1){1.4}}\put(6,2.5){\line(1,-1){1.4}}
\put(11.9,2.6){\line(-1,1){1.35}}\put(11.9,2.4){\line(-1,-1){1.35}}
\put(6,2.5){\line(1,0){5.9}}

\put(10.5,3.865){\line(-1,0){3}} \put(10.5,4.135){\line(-1,0){3}}
\put(10.5,.865){\line(-1,0){3}} \put(10.5,1.135){\line(-1,0){3}}

\put(7.5,4.4){\makebox(0,0){$^1$}}\put(10.5,4.4){\makebox(0,0){$^4$}}\put(7.5,.4){\makebox(0,0){$^2$}}
\put(10.5,.4){\makebox(0,0){$^5$}}\put(5.5,2.5){\makebox(0,0){$^6$}}\put(12.5,2.5){\makebox(0,0){$^3$}}
\put(8.8,4.035){\makebox(0,0){$^\succ$}}\put(9.2,3.765){\makebox(0,0){$^\succ$}}
\put(8.8,1.035){\makebox(0,0){$^\succ$}}\put(9.2,.765){\makebox(0,0){$^\succ$}}

\put(8.8,4.6){\makebox(0,0){$2$}}\put(9.2,3.3){\makebox(0,0){$3$}}
\put(8.8,1.6){\makebox(0,0){$1$}}\put(9.2,.3){\makebox(0,0){$4$}}

\put(7,1.515){\makebox(0,0){$\nwarrow$}}\put(11.505,2){\makebox(0,0){$\swarrow$}}
\put(6.4,1.45){\makebox(0,0){$4$}}\put(11.8,1.6){\makebox(0,0){$1$}}
\end{picture}$$ As in \cite{PBMOG}, orientation and voltage are omitted if an edge gets voltage $0 \in \mathbb{Z}_5$.
}\end{ex}
\medskip

\begin{ex}{\rm
Reye's $(12_4,16_3)$ configuration (cf. e.g. \cite[Figure
2]{BoPi}, \cite[Footnote on P. 140]{HiCV}) is represented by the
$\mathbb{Z}_4$--scheme
$$\left(\begin{smallmatrix}
0&&0&&0&&0\\0,1&&&&2&&3\\&&0,3&&2&&1\\\end{smallmatrix}\right)^{(4)}.$$
}\end{ex}
\medskip

Inverting the approach, $(n_k)$ configurations can be constructed
for whose parameters $n,k$ no instances were known so far:

\begin{ex}{\rm \cite{FLN} The $\mathbb{Z}_7$--scheme
$$
T_{98}^{(7)} \quad = \quad \left(
\begin{smallmatrix}
 0,1,3 & {} & {} & {} & {} & {} & {}&0&0&0&0&0&0&0\\
 {} & 0,1,3 & {} & {} & {} & {} & {}&0&1&2&3&4&5&6\\
 {} & {} & 0,1,3 & {} & {} & {} & {}&0&2&4&6&1&3&5\\
 {} & {} & {} & 0,1,3 & {} & {} & {}&0&3&6&2&5&1&4\\
 {} & {} & {} & {} & 0,1,3 & {} & {}&0&4&1&5&2&6&3\\
 {} & {} & {} & {} & {} & 0,1,3 & {}&0&5&3&1&6&4&2\\
 {} & {} & {} & {} & {} & {} & 0,1,3&0&6&5&4&3&2&1\\
0&0&0&0&0&0&0&2,3,5 & {} & {} & {} &
{} & {} & {}\\0&6&5&4&3&2&1&{} & 2,3,5 & {} & {} & {} & {} & {}\\
0&5&3&1&6&4&2&{} & {} & 2,3,5 & {} & {} & {} &
{}\\0&4&1&5&2&6&3& &{} & {}  & 2,3,5 & {} & {} & {}\\
0&3&6&2&5&1&4&{} & {} & {} & {} & 2,3,5 & {} & {}\\
0&2&4&6&1&3&5&{} & {} & {} & {} & {} & 2,3,5 &
{}\\0&1&2&3&4&5&6 &{} & {} & {} & {} & {} & {} & 2,3,5\\
\end{smallmatrix}
\right)^{(7)}
$$
represents a $(98_{10})$ configuration. Criterion
 \ref{criterion 1} guarantees that the blow--up
$\overline{T_{98}^{(7)}}$ is $J_2$--free: all non--simple full $2
\times 2$ sub--schemes are of type $\left(\smallmatrix 0,1,3 & z\\
-z & 2,3,5\\\endsmallmatrix \right)^{(7)}$  for some $z \in
\mathbb{Z}_7$ and $$ a - z + c + z  \not \equiv 0 \;
(\makebox{mod}\; 7)\;\; \makebox{for all}\; \;a \in
\{0,1,3\}\;\;\makebox{and}\;\; c \in \{2,3,5\}\,,$$ whereas all
simple full $2 \times 2$ sub--schemes arise from the full
multiplication table of $GF(7)$, for details see \cite{FLN}.
}\end{ex}
\medskip

However, not every configuration which has found some
consideration in the literature is polycyclic.

\begin{ex}{\rm
The unique $(9_4,12_3)$ configuration cannot be represented by any
$\mathbb{Z}_3$--scheme. Geometrically, it turns up as the
configuration of the nine points of inflection of a third--order
plane curve without double points in the complex projective plane,
see e.g. \cite[P. 102]{HiCV}. It can also be seen as the affine
plane over $GF(3)$. Its automorphism group has order $432$. A
representation which exhibits a maximum polycyclic
subconfiguration isomorphic to the Pappian $(9_3)$ reads
$$\left(\begin{smallmatrix} 0&&0&&0&&{\bf c_1}\\0&&1&&2&&{\bf c_2}\\0&&2&&1&&{\bf
c_3}\\\end{smallmatrix}\right)^{(3)}$$ where the blow--up
$\overline{\bf c_i}$ of the symbol $\bf c_i$ is the $3 \times 3$
matrix whose entries in the $i^{th}$ column are $1$, and $0$
otherwise (cf. Definition \ref{r-c}).
}\end{ex}
\medskip

\section{Elliptic Semiplanes as Polycyclic Configurations}\label{ellsem}

({\it Desarguesian}) {\it elliptic semiplanes}  are surveyed in
the very last paragraph of Section \ref{verylast}. Let $q = p^\nu$
be a prime power. In \cite{AFLN} it is pointed out that
Desarguesian elliptic semiplanes of types $C$ and $L$ admit
incidence matrices of orders $q^2$ and $q^2 -1 = (q+1)(q-1)$,
respectively, which are $q \times q$ and $(q+1) \times (q+1)$
block matrices with square blocks of orders $q$ and $q-1$. The
blocks are related to certain addition and multiplication tables
of the finite field $GF(q)$. This result has been obtained by
choosing suitable coordinates, which, in turn, depend on the
choice of a suitable labelling for the elements of $GF(q)$. In
general, however, the matrices constructed in \cite{AFLN} cannot
be represented by $\mathbb{Z}_p$--schemes or
$\mathbb{Z}_q$--schemes. In this Section we show how this can be
achieved by fine--tuning the choice of the labelling.

Recurrently we will use the following tool:

\begin{defin} Let $M$ be a matrix of order $(m,n)$  with entries in $GF(q)$.
For each $x \in GF(q)$, we extract its {\bf
position matrix} $P_x$, i.e. the $(0,1)$--matrix of order $(m,n)$
whose entry in position $(i,j)$ is defined by $$ (P_x)_{i,j} \;:=
\; \left\{\begin{array}{cl} 1 & \makebox{if  }x \makebox{ appears as an entry in position } (i,j)\,;
\\ 0 & \makebox{otherwise} \,.
\end{array}\right.$$
\end{defin}

\begin{constr}\label{gfqstar} {\rm
The multiplicative group $(GF(q)^\ast, \cdot)$ is a cyclic group
of order $q-1$, hence one has
$$GF(q)^\ast = \langle y \rangle = \{y, y^2, \ldots,
y^{q-2},y^{q-1}=1\}$$ for a fixed generator $y \in GF(q)^\ast$.
Write down the quotient table of $(GF(q)^\ast, \cdot)$ with
respect to the canonical order $ 1, y, y^2, \ldots, y^{q-2} $ for
the elements:
$$\begin{array}{c|ccccccc|}
:&1&y&y^2&y^3&y^4&\ldots&y^{q-2}\\
\hline1&1&y^{-1}&y^{-2}&y^{-3}&y^{-4}&\ldots&y^{-q+2}\\
y&y&1&y^{-1}&y^{-2}&y^{-3}&\ldots&y^{-q+3}\\
y^2&y^2&y&1&y^{-1}&y^{-2}&\ldots&y^{-q+4}\\
y^3&y^3&y^2&y&1&y^{-1}&\ldots&y^{-q+5}\\
y^4&y^4&y^3&y^2&y&1&\ldots&y^{-q+6}\\
\vdots&\vdots&\vdots&\vdots&\vdots&\vdots&\ddots&\vdots\\
y^{q-2}&y^{q-2}&y^{q-3}&y^{q-4}&y^{q-5}&y^{q-6}&\ldots&1\\\hline
\end{array}$$ Taking into account that $y^{q-1} = 1$ and hence $y^{-q+2} = y$, $y^{-q+3} = y^2$, etc, this
quotient table reveals itself as  a circulant matrix of order
$q-1$. Since an element $x \in GF(q)^\ast$ appears in each row and
column of the quotient table precisely once, its position matrix
$P_x$ is a permutation matrix; in particular, $P_x$ is a circulant
$(0,1)$--matrix, which can be characterised by the only entry $1$
in its first row using the bijection $\iota$; this leads to the
rule
$$P_{y^{-i}} = \overline{\{ i\}} \quad \makebox{for} \; i \in
\mathbb{Z}_{q-1}\,.$$}
\end{constr}

\begin{constr}\label{Lq-1} {\rm Consider the additive group $(GF(q), +)$ and label its
elements, say $x_0, x_1, \ldots, x_{q-1}$, such that $x_0 = 0$.
Write down the difference table of $(GF(q), +)$ with respect to
this labelling. Note that an entry of this difference table is
equal to $0$ if and only if it lies in its main diagonal, whereas
all the other entries are actually elements of $GF(q)^\ast$. Let
$L^{(q-1)} := (\lambda_{ij})_{0 \le i,j \le q}$ be the
$\mathbb{Z}_{q-1}$--scheme of order $q+1$ defined by
$$ \lambda_{i,j}:=  \left\{\begin{array}{ccl} \makebox{blank} &
\makebox{if} & i = j ; \\0 & \makebox{if} &  i = q \;\;
\makebox{or}\; \; j = q\,,\;\makebox{but not both;}\\z &
\makebox{if} & i,j \in
\{0,\ldots,q-1\}\;\makebox{with} \; i \ne j\;\makebox{such that}\;x_i-x_j = y^z.\\
\end{array}\right.$$ Obviously, $L^{(q-1)}$ is a simple $\mathbb{Z}_\mu$--scheme of valency $q$ which
has blank entries on its main diagonal. }\end{constr}
\medskip

\begin{lemma} The $\mathbb{Z}_{q-1}$--scheme $L^{(q-1)}$ is $J_2$--free.
\end{lemma}

\Prf Apply Criterion \ref{criterion 1}: let $\left(\smallmatrix
a&b\\d&c\\\endsmallmatrix\right)^{(q-1)}$ be a full sub--scheme of
$L^{(q-1)}$ and distinguish two cases.

$(i)$ The entries $b$, $c$, and $d$  lie neither in the last
column nor in the last row of $L^{(q-1)}$. Then, by construction,
there exist elements $x_i, x_j, x_g, x_h \in GF(q)$ with $x_i \ne
x_g$ and $x_j \ne x_h$ such that
$$ x_i - x_j = y^a\,,\;\;x_i-x_h = y^b\,,\;\;x_g-x_j =
y^d\,,\;\;x_g-x_h=y^c\;.$$ Then $$a - b + c - d \; \not \equiv \;
0 \;(\makebox{mod } q-1)$$ if and only if
$$1 \ne y^{a-b+c-d} = \frac{y^ay^c}{y^by^d}  = \frac{(x_i -
x_j)(x_g-x_h)}{(x_i-x_h )(x_g-x_j)}=
\frac{x_ix_g-x_ix_h-x_jx_g+x_jx_h}{x_ix_g-x_ix_j-x_hx_g+x_hx_j}\,,
$$ which holds true if and only if $$x_ix_h+x_jx_g \ne x_ix_j+x_hx_g\,,$$
or, equivalently, $$(x_i-x_g)(x_h-x_j) \ne 0\,.$$

$(ii)$ The entries $b$ and $c$  lie in the last column or  $c$ and
$d$ lie in the last row of $L^{(q-1)}$. Then the full $2 \times 2$
sub--scheme reads either   $\left(\smallmatrix
a&0\\d&0\\\endsmallmatrix\right)^{(q-1)}$ or $\left(\smallmatrix
a&b\\0&0\\\endsmallmatrix\right)^{(q-1)}$ and one has $a-b+c-d
\not \equiv 0$ (mod $q-1$) since otherwise either $y^a = y^d$ or
$y^a = y^b$ would appear twice in one and the same column or row
of the difference table $(GF(q),+)$, a contradiction.\qed
\medskip

\begin{prop}
The Desarguesian elliptic semiplane ${\cal S}^L_{q^2-1}$ of type
$L$ derived from $PG(2,q)$  is isomorphic to the $(q-1)$--cyclic
configuration of type $\big((q^2-1)_q \big)$ represented by the
$J_2$-free simple $\mathbb{Z}_{q-1}$--scheme $L^{(q-1)}$ of order
$q+1$ and valency $q$.
\end{prop}

\Prf It is sufficient to check that the construction of
\cite{AFLN} applies also for the above labelling for the elements
of $GF(q)^\ast$. The multiplication table of $(GF(q)^\ast,
\times)$ used in \cite{AFLN}  is actually a quotient table and
matches with the above way of writing it down. \qed
\medskip

\begin{ex}\label{48}{\rm
For later application, construct the $\mathbb{Z}_6$--scheme
$L^{(6)}$ representing an incidence matrix for the Desarguesian
elliptic semiplane ${\cal S}^L_{48}$ on $48$ points. The tables
$$\begin{array}{c|cccccc|}
:&1&3&
2&6&4&5\\
\hline
1&1&5&4&6&2&3\\
3&3&1&5&4&6&2\\
2&2&3&1&5&4&6\\
6&6&2&3&1&5&4\\
4&4&6&2&3&1&5\\
5&5&4&6&2&3&1\\
\hline
\end{array}\qquad\makebox{and}\qquad\begin{array}{c|ccccccc|}
-&0&1&2&
3&4&5&6\\
\hline0&0&6&5&4&3&2&1\\
1&1&0&6&5&4&3&2\\
2&2&1&0&6&5&4&3\\
3&3&2&1&0&6&5&4\\
4&4&3&2&1&0&6&5\\
5&5&4&3&2&1&0&6\\
6&6&5&4&3&2&1&0\\\hline
\end{array}$$
are a quotient table of $GF(7)^\ast = \langle 3 \rangle$ and a
difference table of $GF(7)$ according to Constructions
\ref{gfqstar} and \ref{Lq-1}  , respectively. Then the position
matrices of the elements $1 = 3^0,\; 2 = 3^2,\; 3 = 3^1, \;4 =
3^4,\; 5 = 3^5,$ and $6 = 3^3$ in $GF(7)^\ast$ extracted from the
quotient table read $\overline 0$, $\overline 4$, $\overline 5$,
$\overline 2$, $\overline 1$, and $\overline 3$, respectively, and
the difference table gives rise to the the following
$\mathbb{Z}_6$--scheme:
$$ L^{(6)} \quad = \left(
\begin{smallmatrix}
&3&1&2&5&4&0&0\\
0&&3&1&2&5&4&0\\
4&0&&3&1&2&5&0\\
5&4&0&&3&1&2&0\\
2&5&4&0&&3&1&0\\
1&2&5&4&0&&3&0\\
3&1&2&5&4&0&&0\\
0&0&0&0&0&0&0&\\
\end{smallmatrix}
\right)^{(6)} $$
}\end{ex}

\begin{constr} {\rm Consider the finite field $GF(q)$ as $GF(p)[t]/(f(t))$ for
some irreducible polynomial $f(t) \in GF(p)[t]$ of degree $\nu$.
Then each element in $GF(q)$ can be represented as a polynomial of
degree at most $\nu-1$ with coefficients in $GF(p)$. Label all the
polynomials $\sum_{i=1}^{\nu-1} a_it^i$ with zero constant terms
by $\pi_1 = 0, \pi_2,\ldots,\pi_{p^{\nu-1}}$; they form a subgroup
$S$ of $(GF(q),+)$, which has a copy of $GF(p)$ as direct
complement, namely the constant polynomials. Hence each element in
$GF(q)$ may be written as $\pi_i + z$ for some $\pi_i \in S$ and
$z \in GF(p)$. Choose the canonical order $0,1,\ldots, p-1$ for
the elements of $GF(p)$ and introduce a lexicographic order for
$PG(q)$ by the rule
$$\pi_i + z < \pi_j +w \quad\makebox{if and only if}\quad
\left\{ \begin{array}{l} \makebox{ either}\; i < j \\\makebox{
or}\; \;i = j \;\makebox{and}\; z < w\\\end{array}\right.\,.$$
Write down the difference table of $(GF(q),+)$. Then the block,
say $B_{ij}$, corresponding to minuends  in $\pi_i+GF(p)$ and
subtrahends in $\pi_j+GF(p)$ reads:
$$\begin{array}{c|cccc|}
-&\pi_j&\pi_j+1&\ldots&\pi_j+p-1\\
\hline \pi_i&\pi_i-\pi_j&\pi_i-\pi_j-1&\ldots&\pi_i-\pi_j-p+1\\
\pi_i+1&\pi_i-\pi_j+1&\pi_i-\pi_j&\ldots&+\pi_i-\pi_j-p+2\\
\vdots&\vdots&\vdots&\ddots&\vdots\\
\pi_i+p-1&\pi_i-\pi_j+p-1&\pi_i-\pi_j+p-2&\ldots&\pi_i-\pi_j\\
\hline
\end{array}$$ The block $B_{ij}$ is a circulant matrix, which is
immediately seen by introducing the block
 $$A := \begin{array}{|cccc|}\hline
 0&-1&\ldots&-p+1\\
1&0&\ldots&-p+2\\
\vdots&\vdots&\ddots&\vdots\\
p-1&p-2&\ldots&0\\ \hline
\end{array} \equiv \begin{array}{|cccc|}\hline
 0&p-1&\ldots&1\\
1&0&\ldots&2\\
\vdots&\vdots&\ddots&\vdots\\
p-1&p-2&\ldots&0\\ \hline
\end{array}
$$ (entries taken modulo $p$) and re--writing $B_{ij}$  as
$B_{ij} = \pi_i - \pi_j + A$. With these data, the difference
table becomes a $p^{\nu-1} \times p^{\nu-1}$ block matrix with
circulant blocks of order $p$, namely:
$$\begin{array}{c|cccc|}
-&GF(p)&\pi_1+GF(p)&\ldots&\pi_{p^{\nu-1}}+GF(p)\\
\hline GF(p)&A&-\pi_1+A&\ldots&-\pi_{p^{\nu-1}}+A\\
\pi_1+GF(p)&\pi_1+A&A&\ldots&\pi_1-\pi_{p^{\nu-1}}+A\\
\vdots&\vdots&\vdots&\ddots&\vdots\\
\pi_{p^{\nu-1}}+GF(p)&\pi_{p^{\nu-1}}+A&\pi_{p^{\nu-1}}-\pi_1+A&\ldots&A\\
\hline
\end{array}$$
The block structure reveals the difference table $D_S = (\pi_i -
\pi_j)_{1 \le i,j \le p^{\nu-1}}$ for the subgroup $S$, seen as a
set of representatives for the factor group $GF(q)/GF(p)$ where
$GF(p)$ plays the r\^ ole of the kernel under the epimorphism
$$\epsilon \quad:\quad \left\{
\begin{array}{rcl}(GF(q),+) & \longrightarrow &S\\\sum_{i =
0}^{\nu-1}a_it^i \;+\;(f(t)) &\longmapsto &\sum_{i =
1}^{\nu-1}a_it^i
\\\end{array}\right..$$ We use this fact to construct a
$\mathbb{Z}_p$--scheme $P_{\pi_i + z}$ of order $p^{\nu-1}$ for
each $\pi_i + z \in GF(q)$: extract the position matrix, say
$Q_{\pi_i}$, from the difference table $D_S$ for the group $S$;
then $P_{\pi_i + z}$ is obtained from $Q_{\pi_i}$ by substituting
$\{z\}$ and a blank entry for each entry $1$ and $0$ in
$Q_{\pi_i}$, respectively.}
\end{constr} \medskip

\begin{lemma} The blow--up of the $\mathbb{Z}_p$--scheme $P_{\pi_i + z}$
is the position matrix of the element $\pi_i+z$ extracted from the
above difference table for $GF(q)$. \qed\end{lemma}  \medskip

\begin{constr}
{\rm  Take up the quotient table for $GF(q)^\ast$ from
Construction \ref{gfqstar} and add a new $q^{th}$ row and column
all of whose entries are $0$. Denote the resulting matrix by $G =
(\gamma_{ij})_{0 \le i,j \le q-1}$. Compose a block
$\mathbb{Z}_p$--scheme $C^{(p)} := (\varGamma_{ij})_{0 \le i,j \le
q-1}$ following the rule
$$\varGamma_{ij} \;:=\; P_{\pi+z} \quad \makebox{if and only if}
\quad \gamma_{ij} = \pi + z \in S \oplus GF(p)\,.$$ Seen as a
$\mathbb{Z}_p$--scheme, $C^{(p)}$ is simple and has order $q
p^{\nu-1} = p^{2\nu-1}$ and valency $q$.}
\end{constr}
\medskip

\begin{lemma} The $\mathbb{Z}_{p}$--scheme $C^{(p)}$ is $J_2$--free.
\end{lemma}

\Prf Apply Criterion \ref{criterion 1}: let $\left(\smallmatrix
a&b\\d&c\\\endsmallmatrix\right)^{(p)}$ be a full sub--scheme of
$C^{(p)}$. By construction,  $a,b,c$, and $d$ lie in precisely
four distinct blocks of $C^{(p)}$, say  in $\varGamma_{ij}$, $
\varGamma_{ih}$, $\varGamma_{gh}$, and  $ \varGamma_{gj}$,
respectively, for $i,j,g,h \in \{0,\ldots,q-1\}$ with $i \ne g$
and $j \ne h$. Then there exist elements, say $\pi_a, \pi_b,
\pi_c, \pi_d \in S$, such that $$ \varGamma_{ij} =
P_{\pi_a+a}\,,\quad \varGamma_{ih} = P_{\pi_b+b}\,,
\quad\varGamma_{gh} = P_{\pi_c+c}\,, \quad \makebox{and}
\quad\varGamma_{gj} = P_{\pi_d+d}\,.$$ The upper left element $a$
in $\left(\smallmatrix a&b\\d&c\\\endsmallmatrix\right)^{(p)}$
turns up in $P_{\pi_a+a}$ in local position, say $(i',j')$ for
some $i',j'\in \{1,\ldots,p^{\nu-1}\}$. Analogously, $b, c$, and
$d$ appear in $P_{\pi_b+b}$, $P_{\pi_c+c}$, and  $P_{\pi_d+d}$ in
local positions $(i',h')$, $(g',h')$, and $(g',j')$, respectively,
for some $g',h'\in \{1,\ldots,p^{\nu-1}\}$. This implies that
$\pi_a$, $\pi_b$, $\pi_c$, and $\pi_d$, turn up as entries in the
difference table $D_S$ in positions $(i',j')$, $(i',h')$,
$(g',h')$, and $(g',j')$, respectively.

Hence $$\pi_a = \pi_{i'} - \pi_{j'}\,, \quad\pi_b = \pi_{i'} -
\pi_{h'}\,, \quad\pi_c = \pi_{g'} - \pi_{h'}\,, \quad\pi_d =
\pi_{g'} - \pi_{j'}\,,
$$ and one has $$\pi_a - \pi_b + \pi_c - \pi_d \;= \;\pi_{i'} - \pi_{j'}-
\pi_{i'} + \pi_{h'}+ \pi_{g'} - \pi_{h'}- \pi_{g'} + \pi_{j'}\;
=\; 0\,.$$

Now we distinguish two cases.

$(i)$ The entries $b$, $c$, and $d$  lie neither in the last
column nor in the last row of blocks of $C^{(p)}$, i.e. $i,j,g,h
\in \{0,\ldots,q-2\}$. Hence, by construction,
$$\pi_a+a=\gamma_{ij}= \frac{y^i}{y^j},\; \pi_b+b=
\gamma_{ih}=\frac{y^i}{y^h},
\;\pi_c+c=\gamma_{gh}=\frac{y^g}{y^h},\; \pi_d+d=\gamma_{gj}=
\frac{y^g}{y^j}\,.$$ But then $i \ne g$ and $j \ne h$ imply
$$0 \;\ne \;\frac{(y^i -y^g)(y^h - y^j)}{y^jy^h} \; = \;\frac{y^i}{y^j} - \frac{y^i}{y^h} + \frac{y^g}{y^h} - \frac{y^g}{y^j} \;=$$
$$= \;\pi_a +a - \pi_b -b + \pi_c +c -\pi_d-d \;= \; a-b+c-d\,.$$

$(ii)$ The entries $b=0$ and $c=0$  lie in the last column or $c=$
and $d=0$ lie in the last row of blocks of $C^{(p)}$, i.e. either
$h = q-1$ or $g = q-1$. If $h = q-1$, one has
$$\pi_a+a=\gamma_{ij}= \frac{y^i}{y^j},\quad \pi_b+b=
\pi_c+c=0,\quad \pi_d+d=\gamma_{gj}= \frac{y^g}{y^j}\,,$$ and $i
\ne g$ implies $$0 \;\ne \;\frac{(y^i -y^g)}{y^j} \; =
\;\frac{y^i}{y^j} - \frac{y^g}{y^j} \;=$$
$$= \;\pi_a +a - \pi_b -b + \pi_c +c -\pi_d-d \;= \; a-b+c-d\,.$$
An analogous reasoning works for $g = q-1$. \qed \medskip

\begin{prop} \label{ellsemC}
The Desarguesian elliptic semiplane ${\cal S}^C_{q^2}$ of type $C$
derived from $PG(2,q)$  is isomorphic to the $p$--cyclic
configuration of type $\big((q^2)_q \big)$ represented by the
$J_2$-free simple $\mathbb{Z}_{p}$--scheme $C^{(p)}$ of order
$p^{2\nu-1}$ and valency $q$.
\end{prop}

\Prf Again it is sufficient to check that the construction of
\cite{AFLN} applies also for the above labelling for the elements
of $GF(q)$. \qed
\medskip

\begin{remark}  In {\rm \cite{AFLN}} it has been pointed out that Desarguesian elliptic
semiplanes $\big((q^4-q)_{q^2}\big)$ of type $D$ and Baker's
elliptic semiplane $(45_7)$ of type $B$ do admit representations
by $\mathbb{Z}_{q^2+q+1}$--schemes and a $\mathbb{Z}_3$--scheme,
respectively.
\end{remark}

\section{Regular Graphs of Girth $5$ with Few Vertices}

A $(k,g)$--{\bf cage} is a $k$--regular graph of {\it girth} $g$
with a minimum number of vertices. For a survey on the known
cages, see e.g. in \cite{Wong}. For parameters $k,g$ for which the
$(k,g)$--cage problem is unsolved, some interest has been given to
constructing $k$--regular graphs of girth $g$ with as few vertices
as possible. For $g = 5$, the results given in \cite{AFLN5} have
been outdated by a paper of J\o rgensen's \cite{Joer}, but the
methods based on $\mathbb{Z}_\mu$--schemes presented in
\cite{AFLN5} succeed in tying up with results of \cite{Joer}:

\begin{ex} {\rm The smallest known $9$--regular graph of girth $5$ has $96$ vertices, see \cite{Joer}, Corollary 9.
To construct such a graph, start with the elliptic semiplane
${\cal S}^L_{48}$ on $48$ points. Represent it by the
$\mathbb{Z}_6$--scheme $L^{(6)}$ of order $8$ and valency $7$
constructed in Example \ref{48} and compose the following simple
$\mathbb{Z}_6$--scheme of order $16$ and valency $9$:
$$T_{96}^{(6)}\quad =\quad
\left(
\begin{smallmatrix}
1,5&&&&&&&&\vline&&3&1&2&5&4&0&0\\
&1,5&&&&&&&\vline&0&&3&1&2&5&4&0\\
&&1,5&&&&&&\vline&4&0&&3&1&2&5&0\\
&&&1,5&&&&&\vline&5&4&0&&3&1&2&0\\
&&&&1,5&&&&\vline&2&5&4&0&&3&1&0\\
&&&&&1,5&&&\vline&1&2&5&4&0&&3&0\\
&&&&&&1,5&&\vline&3&1&2&5&4&0&&0\\
&&&&&&&1,5&\vline&0&0&0&0&0&0&0&\\
\hline
&0&2&1&4&5&3&0&\vline&2,4&&&&&&&\\
3&&0&2&1&4&5&0&\vline&&2,4&&&&&&\\
5&3&&0&2&1&4&0&\vline&&&2,4&&&&&\\
4&5&3&&0&2&1&0&\vline&&&&2,4&&&&\\
1&4&5&3&&0&2&0&\vline&&&&&2,4&&&\\
2&1&4&5&3&&0&0&\vline&&&&&&2,4&&\\
0&2&1&4&5&3&&0&\vline&&&&&&&2,4&\\
0&0&0&0&0&0&0&&\vline&&&&&&&&2,4\\
\end{smallmatrix}\right)^{(6)}
$$
The upper right block is a copy of $L^{(6)}$ and the lower left
block is obtained by transposing $L^{(6)}$ and substituting each
entry with its opposite value in $(\mathbb{Z}_6, +)$. Hence both
blocks are $J_2$--free. To check that the whole scheme
$T_{96}^{(6)}$ is $J_2$--free, it is enough to see that each full
$2 \times 2$ subscheme not lying completely in one of these two
blocks is of type
$\left(\smallmatrix 1,5 && z\\
-z & &2,4\\\endsmallmatrix \right)^{(6)}$  for some $z \in
\mathbb{Z}_6$ and $$ a - z + c + z  \not \equiv 0 \;
(\makebox{mod}\; 6)\;\; \makebox{for all}\; \;a \in
\{1,5\}\;\;\makebox{and}\;\; c \in \{2,4\}\,.$$ Since
$T_{96}^{(6)}$ is skew--symmetric and $0$ does not turn up as
entry on its main diagonal, the blow--up $\overline{T_{96}^{(6)}}$
is the adjacency matrix of a $C_4$--free $9$--regular graph $G$ on
$96$ vertices. A short argument shows that $G$  is also
$C_3$--free (cf \cite{AFLN5}, Lemma 2.5). Hence $G$ has girth $
\ge 5$. Equality holds since a $5$-cycle in $G$ is made up by the
vertices corresponding to the $1^{st}, 2^{nd}, 3^{rd}, 91^{st}$,
and $93^{rd}$ rows of $\overline{T_{96}^{(6)}}$.}\end{ex}
\medskip

\begin{remark} The above Example shows  that a $\mathbb{Z}_\mu$--scheme
not only qualifies when major emphasis is laid on an immediate
access to adjacency matrices, but also reveals hidden geometric
structures: consider the {\rm Levi graph} $\Lambda({\cal
S}^L_{48})$, whose adjacency matrix is represented by the
$\mathbb{Z}_6$--scheme $T_{96}^{(6)}$ without its diagonal
entries; this graph is $7$--regular and has girth $6$; then $G$ is
obtained by suitably gluing in $6$-cycles with adjacency matrices
represented by $(1,5)^{(6)}$ and $(2,4)^{(6)}$.
\end{remark}

\begin{ex}{\rm Hoffman-Singleton's celebrated $(7,5)$--cage
\cite{HS60}, say $G_{HS}$, can be obtained in a similar way from
$\Lambda({\cal S}^C_{25})$. In order to construct an adjacency
matrix for $G_{HS}$, we use the representation of $G_{HS}$ due to
Robertson \cite{Ro69}: take five copies $P_0, \ldots, P_4$ of the
pentagram with vertices $0,\ldots,4$ and edges $02,24,41,13,30$,
as well as five copies $Q_0, \ldots, Q_4$ of the pentagon with
vertices $0, \ldots, 4$ and edges $01, 12, 23, 34, 40$. They make
up the $50$ vertices and the first $50$ edges; add further edges
according to the following rule: the vertex $i$ of $P_j$ is joined
to the vertex $l$ of $Q_k$ if, and only if, $$l \equiv i + j k \;
\makebox{(mod}\, 5)\,.$$ Displaying the copies of the pentagrams
and pentagons in the order
$$P_1\,, \ldots, P_4\,,P_0\,,Q_1\,, \ldots, Q_4\,, Q_0 $$ such that the
vertices within each $P_i$ and $Q_j$ maintain the natural  order
$0,1,2,3,4$, the corresponding adjacency matrix turns out to be
the blow--up of the following $\mathbb{Z}_5$--scheme:
$$T_{50}^{(5)}\quad =\quad\left(
\begin{smallmatrix}
2,3&&&&&\vline&1&2&3&4&0\\
&2,3&&&&\vline&2&4&1&3&0\\
&&2,3&&&\vline&3&1&4&2&0\\
&&&2,3&&\vline&4&3&2&1&0\\
&&&&2,3&\vline&0&0&0&0&0\\
\hline
4&3&2&1 &0&\vline&1,4&&&&\\
3&1&4&2&0&\vline&&1,4&&&\\
2&4&1&3&0&\vline&&&1,4&&\\
1&2&3&4&0&\vline&&&&1,4&\\
0&0&0&0&0&\vline&&&&&1,4\\
\end{smallmatrix}\right)^{(5)}$$
The  Levi graph $\Lambda({\cal S}^C_{25})$ has an adjacency matrix
which is represented by the $\mathbb{Z}_5$--scheme $T_{50}^{(5)}$
without its diagonal entries; this graph is $5$--regular and has
girth $6$, and $G_{HS}$ is obtained by suitably gluing in
$5$-cycles with adjacency matrices represented by $(2,3)^{(5)}$
and $(1,4)^{(5)}$.

}\end{ex}
\medskip

\section{Mixed Simple $\mathbb{Z}_\mu$--Schemes}

A configuration $\cal C$ represented by a simple $J_2$--free
$\mathbb{Z}_\mu$--scheme $M^{(\mu)}$ can be partitioned into
$\mu$--sets of pairwise parallel points and lines, say $\Pi_i$ and $\Lambda_j$. A standard
construction in finite geometries applies, namely a kind of projective closure:
new lines $L_i$ and  new points $p_j$ may be added to $\cal C$ such that $L_i$
and $p_j$ are incident with each element in $\Pi_i$ and
$\Lambda_j$, respectively. Eventually, a new point may also be incident with some new line. The following notion renders this
construction compatible with the representation of incidence
matrices as blow--ups of $\mathbb{Z}_\mu$--schemes.

\begin{defin}\label{r-c}
For $ s \ge 1$, we introduce the symbol $\,{\bf r}^s_i$ whose
blow--up is understood to be the $(0,1)$--matrix $\,\overline{{\bf
r}^s_i}$ of order $(s,\mu)$ having entries $1$ in its $i^{th}$ row
and entries $0$ elsewhere. The transpose, denoted by
$\overline{{\bf c}^s_i} := (\overline{{\bf r}^s_i)}^T$, is
interpreted as the blow--up of the symbol $\,{\bf c}^s_i$. Let
$M^{(\mu)} = (z_{ij})$ be a simple $\mathbb{Z}_\mu$--scheme of
order $(m,n)$ with $z_{ij} \in \mathbb{Z}_\mu \cup \{\emptyset\}$.
For permutations $\pi \in S_m$ and $\rho \in S_n$, the scheme
$$M^{(\mu)}_{mix}\quad := \quad \left(
\begin{array}{cccc|c}
z_{11}&z_{12}& \ldots&z_{1n}&{\bf c}^m_{1^\pi}\\
z_{21}&z_{22}& \ldots&z_{2n}&{\bf c}^m_{2^\pi}\\
\vdots&\vdots& \ddots&\vdots&\vdots\\
z_{m1}&z_{m2}& \ldots&z_{mn}&{\bf c}^m_{m^\pi}\\
\hline
{\bf r}^n_{1^\rho}&{\bf r}^n_{2^\rho}& \ldots&{\bf r}^n_{n^\rho}&\bf e\\
\end{array}
\right)^{(\mu)}$$ is called a {\bf mixed}
$\mathbb{Z}_\mu$--scheme, where the blow--up $\overline{\bf e}$ of
the symbol $\bf e$ is a  $(0,1)$--matrix of order $(n,m)$.
\end{defin}

Note that the parameter $\mu$ is not explicitly mentioned in the
symbols ${\bf r}^s_i$ and ${\bf c}^t_j$ since its value coincides
with the parameter $\mu$ of the $\mathbb{Z}_\mu$--scheme under
consideration. In the cases $s = \mu$ and $t = \mu$, we shortly
write ${\bf r}_i$ and ${\bf c}_j$ instead of ${\bf r}^\mu_i$ and
${\bf c}^\mu_j$, respectively. Suitable  matrices $\overline{\bf
e}$ are characterized in the following

\begin{lemma}\label{mixed} Let $M^{(\mu)} = (z_{ij})$ be a simple $J_2$--free
$\mathbb{Z}_\mu$--scheme of order $(m,n)$, where $z_{ij} \in
\mathbb{Z}_\mu \cup \{\emptyset\}$. Then the following are
equivalent
$$\begin{array}{rl}
(i) & \makebox{the blow--up of the mixed
scheme}\;M^{(\mu)}_{mix}\;\makebox{is
still}\;J_2\makebox{--free;}\\
(ii)& \makebox{the blow--up }\overline{\bf e}\makebox{ may have
entry }1 \makebox{ in position }(\pi(j),\rho(i))\makebox{ only if
}z_{ij} = \emptyset.\\
\end{array}
$$
\end{lemma}

\Prf Let the rows and the columns of the blow--up
$\overline{M^{(\mu)}}$ correspond, as usual, to the points and
lines of a configuration $\cal C$. Then the $i^{th}$ column of
$M^{(\mu)}$ gives rise to $\mu$ columns in the blow--up
$\overline{M^{(\mu)}}$. Since $M^{(\mu)}$ is simple, these columns
can be seen as a block matrix made up by just one column of blocks
each of which being either a permutation or a zero matrix of order
$\mu$. Therefore, at most one entry $1$ turns up in each row of
these columns, i.e. any two of the corresponding lines do not have
any point of $\cal C$ in common. Hence these lines make up a
$\mu$-set $\Lambda_i$ of pairwise parallel lines in $\cal C$. An
analogous reasoning holds for any $\mu$-set $\Pi_j$ of points
represented by the $j^{th}$ row of $M^{(\mu)}$. Perform the above
construction and add a new point $p_i$ and a new line $L_j$ such
that $p_i$ and $L_j$ are incident with each element in $\Lambda_i$
and $\Pi_j$, respectively. In terms of incidence matrices, this
means, for each set $\Lambda_i$ and $\Pi_j$, to add a new row and
column to $\overline{M^{(\mu)}}$ which have entries $1$ in
precisely those $\mu$ positions which correspond to the elements
in $\Lambda_i$ and $\Pi_j$, respectively. We distinguish two
cases:

First suppose that $\overline{\bf e}$ is the zero matrix of order
$(n,m)$. This means that no new point lies on any new line. Then
the resulting incidence table is still $J_2$--free. Since this
construction works independently for each row and column of
$M^{(\mu)}$, any permutation $\rho \in S_n$ and $\pi \in S_m$
acting on the indices of the sets $\Lambda_i$ and $\Pi_j$,
respectively, will do. Hence the resulting incidence matrix can be
represented as the blow--up of $M^{(\mu)}_{mix}$ and the
equivalence is clear in this case.

Now suppose that the blow--up $\overline{\bf e}$ has entry $1$ in
position $(i^\rho, j^\pi)$, i.e. the new point $p_{i^\rho}$ is
incident with the new line $L_{j^\pi}$. Then the blow--up
$\overline{M^{(\mu)}_{mix}}$ is $J_2$-free if and only if no line
in $\Lambda_i$ is incident with any point in $\Pi_j$. This, in
turn, is equivalent with $z_{ij} = \emptyset$. Clearly,
$\overline{\bf e}$ is not uniquely determined. \qed
\medskip

\begin{ex}{\rm In Proposition \ref{ellsemC} it has been shown that the Desarguesian
elliptic semiplane ${\cal S}^C_{q^2}$ of type $C$ can be seen as a
$p$--cyclic configuration of type $\big((q^2)_q \big)$. The above
Lemma provides a second representation for ${\cal S}^C_{q^2}$ in
terms of a mixed $\mathbb{Z}_{q-1}$--scheme. Let $C^{(q-1)}$ be
the $\mathbb{Z}_{q-1}$--scheme obtained by deleting the last row
and column in the simple $\mathbb{Z}_{q-1}$--scheme $L^{(q-1)}$
constructed in the proof of Construction \ref{Lq-1}. Since
$C^{(q-1)}$ has blank entries in its main diagonal, Lemma
\ref{mixed} implies that the mixed $\mathbb{Z}_{q-1}$--scheme
$C^{(q-1)}_{mix}$ is $J_2$--free if the blow--up of $\bf e$ is
chosen to be the unit matrix of order $q$. The blow-up
$\overline{C^{(q-1)}_{mix}}$ has valency $q$ and order $q(q-1) + q
= q^2$. }\end{ex}

\begin{remark} The reader will have noticed that the valency of mixed
$\mathbb{Z}_\mu$--schemes has not yet been taken into account.
Obviously, $M^{(\mu)}_{mix}$ has valencies $\mu$ and $\mu+1$ only
if $M^{(\mu)}$ had valencies $\mu - 1$ and  $\mu$, respectively,
and $\overline{\bf e}$ is chosen to be the zero matrix in the
former case and a suitable permutation matrix in the latter case.
On the other hand, {\em partially} mixed $\mathbb{Z}_\mu$--schemes
(i.e. new points and lines are added only for some $\mu$--sets of
pairwise parallel points and lines) can yield
$\mathbb{Z}_\mu$--schemes of valency $k$ even if $M^{(\mu)}$  did
not have a valency. Instances will be discussed in the following
Sections.
\end{remark}

\section{Regular Graphs of Girth $6$ with Few Vertices}

All the known $(k,6)$--cages but one are Levi graphs of finite
projective planes of order $k-1$, the exception being the
$(7,6)$--cage (settled by O'Keefe and Wong \cite{OW}). This cage
revealed itself to be the Levi graph of the elliptic semiplane
$(45_7)$  discovered by Baker some years earlier \cite{Ba78}.
Again, for values $k$ for which the $(k,6)$--cage problem is
unsolved, some interest has been given to finding $k$--regular
graphs of girth $6$ with as few vertices as possible. Levi graphs
of Desarguesian elliptic semiplanes reveal themselves to be good
candidates: for $ k= 11, 13, 16, 19, 23$, and $25$, instances of
smallest known $k$--regular graphs of girth $6$ are $\Lambda({\cal
S}^L_{120})$, $\Lambda({\cal S}^L_{168})$, $\Lambda({\cal
S}^D_{252})$, $\Lambda({\cal S}^L_{360})$, $\Lambda({\cal
S}^L_{528})$, and $\Lambda({\cal S}^D_{620})$, respectively, see
\cite{AFLNA} (cf. also \cite{GaHe}). In \cite{AFLNA}, further
instances have been obtained by deleting an equal number of rows
and columns in $\mathbb{Z}_\mu$--schemes representing Desarguesian
elliptic semiplanes, e.g. a $15$--regular graph on $462$ vertices.
A somewhat more sophisticated and efficient deletion technique in
incidence matrices is due to Balbuena \cite{Balb}, giving rise to
instances of $21$-- and $22$--regular graphs on $964$ and $1008$
vertices, respectively. The methods based on
$\mathbb{Z}_\mu$--schemes presented in \cite{AFLNA} succeed in
tying up with results of \cite{Balb}:

\begin{prop}
For each prime power $q$, there exist $J_2$--free
$(0,1)$--matrices of valency $q-1$ and orders $q^2-q-1$ and
$q^2-q-2$.
\end{prop}

\Prf Consider the simple $\mathbb{Z}_{q-1}$--scheme $L^{(q-1)}$
(see Construction \ref{Lq-1}) and delete two rows and two columns.
In general, this yields a simple $\mathbb{Z}_{q-1}$--scheme of
order $q-1$, which has $q-i$ blank entries, with $i = 1,2,3$. For
the last two cases, we can  choose the following minors
$M^{(q-1)}$ and $N^{(q-1)}$ of $L^{(q-1)}$, which are respectively
obtained by deleting the first two rows as well as
$$\begin{array}{rl}
\makebox{the first and the last columns} & \makebox{if}\; \;i = 2\,,\\
 \makebox{the last two columns} & \makebox{if}\;\; i = 3\,.\\
\end{array}
$$  Embed the $\mathbb{Z}_{q-1}$--schemes $M^{(q-1)}$ and
$N^{(q-1)}$ into the  mixed schemes
$$
M_{mix}^{(q-1)}:= \left(
\begin{array}{c|ccccc|c}
z_{11}&\emptyset&z_{13}&z_{14}&\ldots&z_{1,q-1}&{\bf c}^{q-2}_1\\
z_{21}&z_{22}&\emptyset&z_{24}&\ldots&z_{2,q-1}&{\bf c}^{q-2}_2\\
\vdots&\vdots&\ddots&\ddots&\ddots&\vdots&\vdots\\
z_{q-3,1}&z_{q-3,2}&z_{q-3,3}&\ddots&\emptyset&z_{q-3,q-1}&{\bf c}^{q-2}_{q-3}\\
z_{q-2,1}&z_{q-2,2}&z_{q-2,3}&\ldots&z_{q-2,q-2}&\emptyset&{\bf
c}^{q-2}_{q-2}\\\hline
z_{q-1,1}&z_{q-1,2}&z_{q-1,3}&z_{q-1,4}&\ldots&z_{q-1,q-1}&\emptyset\\\hline
\emptyset&{\bf r}^{q-2}_1&{\bf r}^{q-2}_2&{\bf r}^{q-2}_3&\ldots&{\bf r}^{q-2}_{q-2}&\emptyset\\
\end{array}
\right)^{(q-1)}$$ and $N_{mix}^{(q-1)}:=$
$$
= \left(
\begin{array}{cc|ccccc|c}
z_{11}&z_{12}&\emptyset&z_{14}&z_{15}&\ldots&z_{1,q-1}&{\bf c}^{q-3}_1\\
z_{21}&z_{22}&z_{23}&\emptyset&z_{25}&\ldots&z_{2,q-1}&{\bf c}^{q-3}_2\\
\vdots&\vdots&\vdots&\ddots&\ddots&\ddots&\vdots&\vdots\\
z_{q-4,1}&z_{q-4,2}&z_{q-4,3}&z_{q-4,4}&\ddots&\emptyset&z_{q-4,q-1}&{\bf c}^{q-3}_{q-4}\\
z_{q-3,1}&z_{q-3,2}&z_{q-3,3}&z_{q-3,4}&\ldots&z_{q-3,q-2}&\emptyset&{\bf
c}^{q-3}_{q-3}\\\hline
z_{q-2,1}&z_{q-2,2}&z_{q-2,3}&z_{q-2,4}&z_{q-2,5}&\ldots&z_{q-2,q-1}&\emptyset\\
z_{q-1,1}&z_{q-1,2}&z_{q-1,3}&z_{q-1,4}&z_{q-1,5}&\ldots&z_{q-1,q-1}&\emptyset\\\hline
\emptyset&\emptyset&{\bf r}^{q-3}_1&{\bf r}^{q-3}_2&{\bf r}^{q-3}_3&\ldots&{\bf r}^{q-3}_{q-3}&\emptyset\\
\end{array}
\right)^{(q-1)}$$ The valency of both $M_{mix}^{(q-1)}$ and
$N_{mix}^{(q-1)}$ is $q-1$ and their orders are $$(q-1)(q-1) + q -
i = q^2 - q -i +1 $$ for $i = 2$ and $i=3$, respectively. Then
their blow--ups will do.\qed \medskip

\section{Kr\v cadinac's Configuration of Type $34_6$}

In this Section we present a construction yielding four
configurations of type $30_5$, which will be used to obtain Kr\v
cadinac's configuration of type $(34_6)$ (cf. \cite{krc}) and four
new configurations of type $(35_6)$. The computer results have
been obtained by using the software \cite{Kocay}.

\begin{constr} {\rm Start with the elliptic semiplane  ${\cal
S}^L_{15}$ and represent it by the $\mathbb{Z}_3$--scheme
$L^{(3)}$ of order $5$ and valency $4$, see Construction \ref{Lq-1}.
Compose the following simple $\mathbb{Z}_3$--scheme of order $10$
and valency $5$
$$T =
\left(
\begin{smallmatrix}
\alpha_1&&&&&\vline&&0&1&2&0\\
&\alpha_2&&&&\vline&0&&2&1&0\\
&&\alpha_3&&&\vline&1&2&&0&0\\
&&&\alpha_4&&\vline&2&1&0&&0\\
&&&&\alpha_5&\vline&0&0&0&0&\\\hline
&0&2&1&0&\vline&\beta_1&&&&\\
0&&1&2&0&\vline&&\beta_2&&&\\
2&1&&0&0&\vline&&&\beta_3&&\\
1&2&0&&0&\vline&&&&\beta_4&\\
0&0&0&0&&\vline&&&&&\beta_5\\
\end{smallmatrix}\right)^{(3)}
$$
for suitable $\alpha_i, \beta_i \in \mathbb{Z}_3$. The upper right
block is a copy of $L^{(3)}$ and the lower left block is obtained
by transposing $L^{(3)}$ and substituting each entry by its
opposite element in $(\mathbb{Z}_3,+)$. }\end{constr}

\begin{lemma}  The $\mathbb{Z}_3$--schemes obtained for
$$\begin{array}{lcc}
T_{360}&:&(\alpha_1,\ldots,\alpha_5) = (1,1,1,1,1), \;(\beta_1,\ldots,\beta_5) = (1,1,1,1,1),\\
T_{72}&:&(\alpha_1,\ldots,\alpha_5) = (1,1,1,1,0), \;(\beta_1,\ldots,\beta_5) = (1,1,1,1,0),\\
T_{36}&:&(\alpha_1,\ldots,\alpha_5) = (1,1,1,1,1), \;(\beta_1,\ldots,\beta_5) = (1,1,1,1,0),\\
T_{18}&:&(\alpha_1,\ldots,\alpha_5) = (1,1,1,1,1), \;(\beta_1,\ldots,\beta_5) = (1,1,1,0,0)\\
\end{array}$$
represent four pairwise non-isomorphic configurations ${\cal
T}_{360}$, ${\cal T}_{72}$,${\cal T}_{36}$, and ${\cal T}_{18}$ of
type $(30_5)$, whose automorphism groups have orders $360$, $72$,
$36$, and $18$, respectively. \qed
\end{lemma}

\Prf Apply Criterion \ref{criterion 1} to the
$\mathbb{Z}_3$--scheme $T$: all full $2 \times 2$ sub--schemes are
of type $\left(
\begin{smallmatrix}
\alpha_i&\lambda_{ij}\\
-\lambda_{ij}&\beta_j\\
\end{smallmatrix}\right)^{(3)}$, for $i,j \in \{1, \ldots, 5\}$ with $i \ne j$.
Thus $T$ meets the condition of the criterion if and only if
$$(\ast)\quad \alpha_i + \beta_j \not \equiv 0 \;(\makebox{mod} \;3)
\quad\makebox{for all}\; i,j = 1, \ldots,5\;\makebox{with} \; i
\ne j\,.$$ There are a lot of solutions for $(\ast)$. A computer
search, however,  reveals that they lead to only four pairwise
non--isomorphic configurations. We can choose the solutions
indicated above.
  \qed \medskip

\medskip

\begin{constr}{\rm Let $\cal T$ stand for one of the four configurations ${\cal
T}_{360}$, ${\cal T}_{72}$, ${\cal T}_{36}$, or ${\cal T}_{18}$ of
type $(30_5)$.  Rearrange both the rows and columns of the
$\mathbb{Z}_3$--scheme $T$ following the order $1, 6, 2, 7, 3, 8,
4, 9, 5, 10$, to obtain an equivalent variant, namely
$$V(T) =\left(
\begin{smallmatrix}
1&&|&&0&|&&1&|&&2&|&&0\\
&1&|&0&&|&2&&|&1&&|&0&\\\hline
&0&|&1&&|&&2&|&&1&|&&0\\
0&&|&&1&|&1&&|&2&&|&0&\\\hline
&1&|&&2&|&1&&|&&0&|&&0\\
2&&|&1&&|&&1&|&0&&|&0&\\\hline
&2&|&&1&|&&0&|&1&&|&&0\\
1&&|&2&&|&0&&|&&\beta_4&|&0&\\\hline
&0&|&&0&|&&0&|&&0&|&\alpha_5&\\
0&&|&0&&|&0&&|&0&&|&&\beta_5\\
\end{smallmatrix}\right)^{(3)}.$$
Note that, for $j = 1,3,5,7,9$, the $j^{th}$ and $j+1^{st}$ rows
(and columns) of the scheme $V(T)$ are {\bf non--overlapping},
i.e. the entries in one and the same position of the $j^{th}$ and
$j+1^{st}$ rows (and columns) are always one element of
$\mathbb{Z}_3$ and one blank entry. Hence, in the blow--up
$\overline{V(T)}$ of $V(T)$, the rows (columns) labelled by
$$(\S) \quad \quad 3(j-1)+1, \;\;3(j-1)+2, \;\;3j, \;\;3j+1, \;\;3j+2, \;\;3(j+1)$$
correspond to $6$ pairwise parallel points (lines) of $\cal T$.
Denote the sets of these six points  and lines by $\Pi_l$ and
$\Lambda_l$, respectively, where $l := \frac12 (j+1)$. The
families $\{\Pi_l\}_{l = 1,\ldots,5}$ and $ \{\Lambda_l \}_{l =
1,\ldots,5}$ partition the sets of all points and lines in  ${\cal
T}$. A computer evaluation reveals the following} \end{constr}

\begin{lemma} \label{aut} The families $\{\Pi_l\}_{l
= 1,\ldots,5}$ and $ \{\Lambda_l \}_{l = 1,\ldots,5}$
 are invariant under all automorphisms of ${\cal T}$. \qed
\end{lemma}

\medskip

\begin{constr}{\rm  $\;$ Now let $\cal T$ stand for one of the three
configurations ${\cal T}_{360}$, ${\cal T}_{72}$, and ${\cal
T}_{36}$ of type $(30_5)$, represented by the schemes $V(T)$
obtained by Construction $2$. For $l = 1, \ldots, 4$, add a new
\lq \lq improper" line and point for each set $\Pi_l$ and
$\Lambda_l$. Equivalently, add four new rows and columns to the
blow--up $\overline{V(T)}$ which, for $j = 1,3,5,7$, have entries
$1$ in positions $(\S)$ and entries $0$ else. Simultaneously,
substitute the $2 \times 2$ sub--scheme $\left(
\begin{smallmatrix}
\alpha_5&\\
&\beta_5\\
\end{smallmatrix}\right)^{(3)}$ of $V(T)$ by $ \left(
\begin{smallmatrix}
\alpha_5,\eta&\\
&\beta_5,\zeta\\
\end{smallmatrix}\right)^{(3)}$ for some $\eta, \zeta \in
\mathbb{Z}_3$. The mixed scheme
$$V(T)' =\left(
\begin{smallmatrix}
1&&|&&0&|&&1&|&&2&|&&0&|&{\bf c}^4_1\\
&1&|&0&&|&2&&|&1&&|&0&&|&{\bf c}^4_1\\\hline
&0&|&1&&|&&2&|&&1&|&&0&|&{\bf c}^4_2\\
0&&|&&1&|&1&&|&2&&|&0&&|&{\bf c}^4_2\\\hline
&1&|&&2&|&1&&|&&0&|&&0&|&{\bf c}^4_3\\
2&&|&1&&|&&1&|&0&&|&0&&|&{\bf c}^4_3\\\hline
&2&|&&1&|&&0&|&1&&|&&0&|&{\bf c}^4_4\\
1&&|&2&&|&0&&|&&\beta_4&|&0&&|&{\bf c}^4_4\\\hline
&0&|&&0&|&&0&|&&0&|&\alpha_5, \eta&&|&\\
0&&|&0&&|&0&&|&0&&|&&\beta_5, \zeta&|&\\\hline
{\bf r}^4_1&{\bf r}^4_1&|&{\bf r}^4_2&{\bf r}^4_2&|&{\bf r}^4_3&{\bf r}^4_3&|&{\bf r}^4_4&{\bf r}^4_4&|&&&|&\\
\end{smallmatrix}\right)^{(3)}$$
suitably represents the result of these modifications. Note that
$V(T)'$ has valency $6$.}\end{constr}

\begin{lemma} The $\mathbb{Z}_3$--schemes $V(T_{360})'$, $V(T_{72})'$,
and $V(T_{36})'$ turn out to be $J_2$--free for just one pair
$(\eta, \zeta)$ each, namely $(0,0)$, $(1,1)$, and $(0,1)$,
respectively. All three solutions lead to one and the same mixed
scheme with $\{\alpha_5,\eta\} = \{\beta_5, \zeta\} = \{0,1\}$,
whose blow--up represents Kr\v cadinac's configuration of type
$(34_6)$ {\rm \cite{krc}}. Its automorphism group has order $72$.
\end{lemma}

\Prf A straightforward verification shows that the
$\mathbb{Z}_3$--scheme $V(T)'$ is $J_2$--free. The isomorphism
with Kr\v cadinac's configuration and the order of its
automorphism group have been obtained by computer.\qed

\begin{remark} Let $\cal T$ stand for one of the four
configurations ${\cal T}_{360}$, ${\cal T}_{72}$,  ${\cal
T}_{36}$, or ${\cal T}_{18}$. Alternatively, we can also add five
new \lq \lq improper" lines and points for the families
$\{\Pi_l\}_{l = 1,\ldots,5}$ and $ \{\Lambda_l \}_{l =
1,\ldots,5}$ in $\cal T$, respectively. This leads to four
configurations of type $(35_6)$, represented by the mixed
$\mathbb{Z}_3$--schemes
$$\left(\begin{smallmatrix}
1&&|&&0&|&&1&|&&2&|&&0&|&{\bf c}^5_1\\
&1&|&0&&|&2&&|&1&&|&0&&|&{\bf c}^5_1\\\hline
&0&|&1&&|&&2&|&&1&|&&0&|&{\bf c}^5_2\\
0&&|&&1&|&1&&|&2&&|&0&&|&{\bf c}^5_2\\\hline
&1&|&&2&|&1&&|&&0&|&&0&|&{\bf c}^5_3\\
2&&|&1&&|&&1&|&0&&|&0&&|&{\bf c}^5_3\\\hline
&2&|&&1&|&&0&|&1&&|&&0&|&{\bf c}^5_4\\
1&&|&2&&|&0&&|&&\beta_4&|&0&&|&{\bf c}^5_4\\\hline
&0&|&&0&|&&0&|&&0&|&\alpha_5&&|&{\bf c}^5_5\\
0&&|&0&&|&0&&|&0&&|&&\beta_5&|&{\bf c}^5_5\\\hline
{\bf r}^5_1&{\bf r}^5_1&|&{\bf r}^5_2&{\bf r}^5_2&|&{\bf r}^5_3&{\bf r}^5_3&|&{\bf r}^5_4&{\bf r}^5_4&|&{\bf r}^5_5&{\bf r}^5_5&\\
\end{smallmatrix}\right)^{(3)}$$
whose automorphism groups have still orders $360$, $72$, $36$, and
$18$ (cf. {\rm Lemma \ref{aut}}). These configurations are new.
Thus far, three configurations of type $(35_6)$ have been
exhibited in the literature: In {\rm \cite{Gropp90}} and {\rm
\cite{MPW}}, cyclic configurations are presented in terms of
deficient cyclic difference sets, namely
$${\cal
C}_G \;:\; \{0,1,8,11,13,17\}^{(35)}\quad\makebox{and}\quad{\cal
C}_{MPW}\;:\;\{0,1,3,7,12,20\}^{(35)} \,,$$ respectively, whereas
 in {\rm \cite{FLN}} there is mentioned a
configuration ${\cal C}_{FLN}$  represented by
the following  $\mathbb{Z}_7$--scheme:
$$
 \left(
\begin{smallmatrix}
0,1 & 6 & 2 & 2 & 6 \\
6 & 0,1 & 6 & 2 & 2 \\
2 & 6 & 0,1 & 6 & 2 \\
2 & 2 & 6 & 0,1 & 6 \\
6 & 2 & 2 & 6  & 0,1
\end{smallmatrix}
\right)^{(7)}
$$

A computer check reveals that ${\cal C}_G$ is isomorphic to ${\cal
C}_{MPW}$; its automorphism group has order $35$, whereas ${\cal
C}_{FLN}$ has an automorphism group of order $140$. It is cyclic
as well and isomorphic to the configuration given by the deficient
difference set $\{0,1,8,12,14,17\}^{(35)}$. A computer search
confirms that there are no further cyclic configurations of type
$35_6$.

\end{remark}

\section{Appendix: $(0,1)$-Matrices, Graphs, and Configurations}\label{verylast}

A {\it circulant} matrix is a square matrix where each row vector
is shifted one element to the right relative to the preceding row
vector. Hence a circulant $(0,1)$--matrix is uniquely determined
by the positions of the entries $1$ in its first row. The
transpose of a matrix $A$ is denoted by $A^T$.

Graph theoretic notations come from \cite{Biggs}. We distinguish
{\it graphs} from {\it general graphs}, the former having neither
loops nor multiple edges. All (general) graphs are supposed to be
finite and connected (if not otherwise stated).

Let $K$ be a general graph all of whose edges have been given plus
and minus directions. A {\it cyclic voltage graph} is the pair
$(K,\alpha)$ where $\alpha$ is a function from the $+$ directed
edges of $K$ into the cyclic group $\mathbb{Z}_\mu$, called a {\it
cyclic voltage assignment}. For slightly different and more
general definitions, cf. e.g. \cite{Gross, GT1,GT2, Pisanski,
PBMOG}. The {\it derived graph} $K^\alpha$, also referred to as
the {\it lift} of K in $\mathbb{Z}_\mu$ via $\alpha$ (cf. e.g.
\cite{Exoo}) or the ({\it regular}) {\it covering graph} (cf. e.g.
\cite{Pisanski,White}), is the (not necessarily connected) general
graph whose vertex and edges sets are $VK \times \mathbb{Z}_\mu$
and $EK \times \mathbb{Z}_\mu$ and in which $(v,a)$ and $(w,b)$
are incident with $(e,a)$ if $EK$ contains an edge $e$ whose $+$
direction runs from $v$ to $w$ and $a + \alpha(e) = b$.  Note that
\lq\lq regular" has a topological meaning (cf. e.g.
\cite{GT1,GT2}). The {\it natural projection} $\pi : K^\alpha
\longrightarrow K$ is defined by the rules $(v,a)^\pi = v$ and
$(e,a)^\pi = e$.

\begin{lemma} The lift of K in $\mathbb{Z}_\mu$ via
$\alpha$ is a graph if loops of $K$ don't have image $0 \in
\mathbb{Z}_\mu$ and multiple edges do have distinct images under
the cyclic voltage assignment.
\end{lemma}

\Prf Let $e$ be a $v$--based loop in $K$ with voltage $a \in
\mathbb{Z}_\mu \backslash \{0\}$. If $a$ has order $\nu$, then the
loop gives rise to $\frac\mu\nu$ cycles of length $\nu$ in
$K^\alpha$, namely
$$(v,c), (v,c+a), (v,c+2a), \ldots, (v,c+(\nu-1)a)$$ for $c =
0,\ldots, \frac\mu\nu -1$. Let $e,f \in EK$ be a double edge in
$K$, both running from $v$ to $w$, with voltages $a,b$,
respectively. This leads to $2\mu$ distinct edges in $K^\alpha$,
no two of which incident with the same pair of vertices,
namely
$$(v,c) | (e,a) | (w,c+a) \quad \makebox{and} \quad (v,c) |
(f,b) | (w,c+b)$$ for $c \in \mathbb{Z}_\mu$. \qed \medskip

In the light of this Lemma, we may call a cyclic voltage
assignment $\alpha: K \longrightarrow \mathbb{Z}_\mu$ {\it
admissible} if loops of $K$ don't have image $0 \in
\mathbb{Z}_\mu$ and multiple edges do have distinct images.

Suppose that $\varGamma$ is a graph whose vertex set $V\varGamma$
is the set $\{v_1,\ldots,v_n\}$, and consider the edge set
$E\varGamma$ as a set of unordered pairs of elements in
$V\varGamma$: then the {\it adjacency matrix} of $\varGamma$ is
the $n \times n$ matrix $A = A(\varGamma)$ whose entries $a_{ij}$
are given by $a_{ij} := 1$ if $\{v_i, v_j\} \in E\varGamma$, and
$a_{ij} := 0$ otherwise. $A$ is a symmetric matrix with entries
$0$ on the main diagonal. The rows and columns of $A$ correspond
to an arbitrary labelling of the vertices of $\varGamma$. A
permutation $\pi$ of $V\varGamma$ can be represented by a {\it
permutation matrix} $P_\pi = (p_{ij})$, where $p_{ij} = 1$ if $v_i
= v_j^\pi$, and $p_{ij} = 0$ otherwise. Then $P_\pi^{-1}AP_\pi$
becomes the adjacency matrix of $\varGamma$ with respect to this
re--labelling. Thus we focus primarily on the equivalence class
$\cal A$ of $(0,1)$--matrices represented by $A$ under the
equivalence relation
$$A_1 \cong A_2 \quad \makebox{if} \quad  A_2 = P_\pi^{-1}A_1P_\pi \;\;
\makebox{for some permutation matrix}\; P_\pi \; \makebox{with}\;
\pi \in S_n
$$ on the set of symmetric $(0,1)$--matrices with zero diagonal.

A graph is called $k$--{\it regular} if every vertex is adjacent
to $k$ distinct vertices.  A graph is {\it bipartite} if its vertex set
can be partitioned into two parts $V_1$ and $V_2$ such that each
edge has one vertex in $V_1$ and one vertex in $V_2$. If we label
the vertices in such a way that those in $V_1$ come first, then
the adjacency matrix of a bipartite graph takes the form $$A =
\begin{pmatrix} 0 & B\\B^T & 0\end{pmatrix}\,.$$ A bipartite graph
is $(k,l)$--{\it semiregular} if the vertex in $V_1$ and $V_2$ are
adjacent to $k$ and $l$ vertices, respectively.

An adjacency matrix for the  {\it cycle (graph)} ${\cal C}_n$ is
the circulant matrix with first row $[0,1,0,\ldots,0,1]$. The {\it
girth} of a graph $\varGamma$ is the length $g$ of a shortest
cycle ${\cal C}_g$ which can be embedded into $\varGamma$.

\begin{lemma}\label{girth} {\rm [Folk]} Let $\varGamma$
be a bipartite graph. Then the following are equivalent:

$$\begin{array}{rl}
$(i)$ & \varGamma \; \makebox{has girth}\; \ge 6;\\
$(ii)$ & \varGamma \; \makebox{is}\;\; {\cal C}_4\makebox{--free};\\
$(iii)$ & \makebox{the adjacency matrix}\; A(\varGamma) \;
\makebox{is}\;\; J_2\makebox{--free}.\end{array}
$$ \qed
\end{lemma}

Out of the many ways to introduce configurations, we prefer Levi's
definition \cite{Levi}, which best suits our approach to Graph
Theory via $(0,1)$--matrices. An {\it incidence table} or {\it
incidence matrix} $C$ is a $J_2$--free $(0,1)$--matrix; usually
some {\it regularity} is requested:  $C$ is of {\it type}
$(m_k,n_l)$ if $C$ has order $(m,n)$ and
if the sums of all
entries in the rows and columns have constant values $k$ and $l$,
respectively.
The meaning
of {\it points, lines, incidences } etc. are based on the usual
interpretation of an incidence table. A {\it schematic
configuration} $(m_k,n_l)$ is an equivalence class $\cal C$ of
incidence tables of type $(m_k,n_l)$ under the equivalence
relation $$C_1 \cong C_2 \quad \makebox{if} \quad  C_2 = PC_1Q
\;\; \makebox{for permutation matrices}\;  P \; \makebox{and} \; Q
\,.$$ Other names are {\it combinatorial configuration} or simply
{\it configuration}, not to be confused with a {\it geometric
configuration} made up by points and lines of the Euclidean plane.
If $m = n$ (and hence $k = l$), the symbol $(n_k,n_k)$ will be
shortened to $(n_k)$. In the literature, such configurations are
called {\it symmetric}. We avoid this term, since \lq \lq
symmetric" configurations need not admit symmetric incidence
tables.

With each $(m_k,n_l)$ configuration $\cal C$ one associates its
{\it Levi graph} $\Lambda({\cal C})$, see \cite{Cox}: it is the
bipartite graph whose vertices  are the points and lines of $\cal
C$; two vertices of $\Lambda({\cal C})$ are adjacent if and only
if they make up an incident point--line pair in $\cal C$. If $\cal
C$ is represented by the incidence table $C$, then $A :=
\left(\begin{smallmatrix} 0 & C\\C^t & 0\end{smallmatrix}\right)$
is an adjacency matrix for $\Lambda({\cal C})$. Lemma \ref{girth}
implies that Levi graphs of configurations have girth $\ge 6$.

For each $n \in \mathbb{N}$ and $1 \le k \le \frac12 +
\sqrt{n-\frac34}$, a subset $D = \{s_0, \ldots, s_{k - 1} \}
\subseteq \mathbb{Z}_n$ is called a {\it deficient cyclic
difference set}, denoted by $\{s_0, \ldots, s_{k - 1} \}^{(n)}$,
if the $k^2 - k$ differences $s_i - s_j$ (mod $n$)
 are distinct in pairs for $i,j = 0, \ldots, k - 1$ with $i \ne
j$, see e.g. \cite{Funk06,MPW}. The {\it deficiency} $d := n - k^2
+ k - 1$ counts how many elements in $\mathbb{Z}_n^\ast$ are not
covered by any such difference.

A configuration $(n_k)$ is called {\it cyclic} if its points can
be labelled by the elements of $\mathbb{Z}_n$ such that its lines
are given by a {\it base--line}, i.e. a set $\{z_0, \ldots,
z_{k-1}\}$ of $k$ distinct points,  and all its {\it shifts}
$\{z_0+c, \ldots, z_{k-1}+c\}$, numbers taken modulo $n$, for $c =
1, \ldots, n-1$.

\begin{lemma}\label{cds} {\rm \cite{Gropp90,Lipman}} A subset $D
\subseteq \mathbb{Z}_n$ of cardinality $k$ is the base line of
some cyclic configuration $(n_k)$ if and only if $D$ is a
deficient cyclic difference set. \qed
\end{lemma}

A finite {\it elliptic semiplane of order} $k-1$  is an $(n_k)$
configuration satisfying the following axiom of parallels: given a
non-incident point line pair $(p_1,L_1)$, there exists at most one
line $L_2$ incident with $p_1$ and {\it parallel} to $L_1$ (i.e.
there is no point incident with both $L_1$ and $L_2$) and at most
one point $p_2$ incident with $L_1$ and {\it parallel} to $p_1$
(i.e. there is no line incident with both $p_1$ and $p_2$), for
details, see e.g. \cite{Dembowski}.

For a survey on the known examples the following notion is useful:
a {\it Baer subset} of a finite projective plane $\cal P$ is
either a Baer subplane $\cal B$ or, for a distinguished
point--line pair $(p_0,L_0)$, the union ${\cal B}(p_0,L_0)$ of all
lines and points incident with $p_0$ and $L_0$, respectively.
Trivial examples of elliptic semiplanes are finite projective
planes of order $n$, which are $\big((n^2+n+1)_{n+1}\big)$
configurations. Instances (of {\it type $L$}) are obtained from
finite projective planes of order $n$ by deleting a Baer subset
${\cal B}(p_0,L_0)$ where $(p_0,L_0)$ is a distinguished
non--incident point--line pair. The resulting structures are
$\big((n^2-1)_n \big)$ configurations. Similarly, instances (of
{\it type $C$}) are obtained from finite projective planes of
order $n$ by deleting a Baer subset ${\cal B}(p_1,L_1)$ with
$(p_1,L_1)$ incident, yielding $\big((n^2)_n \big)$
configurations. Complements ${\cal P} \backslash {\cal B}$ of Baer
subplanes $\cal B$ make up a third series of instances (of {\it
type $D$}), furnishing $\big((n^4-n)_{n^2} \big)$ configurations.
A sporadic example is the elliptic semiplane $(45_7)$ found by
Baker \cite{Ba78}. Elliptic semiplanes of types $C,D$, and $L$ are
said to be {\it Desarguesian} and denoted by ${\cal S}^C$, ${\cal
S}^D$, and ${\cal S}^L$, respectively, if they are derived from
$PG(2,q)$.

\end{document}